\documentclass[12pt]{amsart}
\usepackage{latexsym,amsmath,amsfonts,amsthm,amssymb}
\usepackage{version}
\begin{document}

\newcommand{\ds}{\displaystyle} 

\newcommand{\R}{{\mathbb R}}
\newcommand{\N}{{\mathbb N}}
\newcommand{\C}{{\mathbb C}}
\newcommand{\Q}{{\mathbb Q}}
\newcommand{\Z}{{\mathbb Z}}
\newcommand{\Rn}{{\mathbb R}^n}
\newcommand{\M}{{ \mathcal H}}

\newcommand{\1}{{{\mathbf 1}}}

\newcommand{\SX}{{ l^1S}}
\renewcommand{\H}{\dot{H}}
\newcommand{\tS}{\tilde{S}}
\newcommand{\tT}{\tilde{T}}
\newcommand{\tW}{\tilde{W}}
\newcommand{\tX}{\tilde{X}}
\newcommand{\tK}{\tilde{K}}

\newcommand{\tx}{\tilde{x}}
\newcommand{\txi}{\tilde{\xi}}

\newcommand{\tA}{\tilde{A}}
\newcommand{\tI}{{\tilde I}}
\newcommand{\tP}{\tilde{P}}

\newcommand{\la}{{\langle}}
\newcommand{\ra}{{\rangle}}

\newcommand{\supp}{{\text{supp }}}

\newcommand{\bc}{\begin{2}\begin{com}}
\newcommand{\ec}{\end{com}\end{2}}

\newcommand{\cd}{{\, \cdot\,}}

\renewcommand{\varepsilon}{\epsilon}
\renewcommand{\div}{{\text{div}\,}}

\newtheorem{theorem}{Theorem}
\newtheorem{lemma}[theorem]{Lemma}
\newtheorem{proposition}[theorem]{Proposition}
\newtheorem{corollary}[theorem]{Corollary}
\newtheorem{definition}[theorem]{Definition}
\newtheorem{remark}[theorem]{Remark}
\newtheorem{assumption}[theorem]{Assumption}

\newcommand{\bb}{\beta}
\newcommand{\e}{\epsilon}
\newcommand{\oo}{\omega}
\newcommand{\dd}{\delta}
\renewcommand{\d}{{\partial}}
\renewcommand{\l}{\lambda}
\newcommand{\p}{{\partial}}

\newenvironment{com}{\begin{quotation}{\leftmargin .25in\rightmargin .25in}\sffamily \footnotesize $\clubsuit$}
               {$\spadesuit$\end{quotation}\par\bigskip}
\newenvironment{com2}{\sffamily\footnotesize $\clubsuit$ }{ $\spadesuit$}

\bibliographystyle{plain}

\title{Global  parametrices and dispersive estimates for 
 variable coefficient wave equations}

\author[J. Metcalfe]
{Jason Metcalfe}

\address{Department of Mathematics, University of North Carolina,
  Chapel Hill, NC  27599-3250, USA}

\author[D. Tataru]
{Daniel Tataru}

\address{Department of Mathematics, University of California \\
Berkeley, CA 94720-3840, USA}

\thanks{ 
The work of the first author was supported in part by an NSF
postdoctoral fellowship and NSF grant DMS0800678, and that of the second author by
NSF grants DMS0354539 and DMS0301122.}

\baselineskip 18pt

\begin{abstract}
  In this article we consider variable coefficient time dependent
  wave equations  in $\R \times \R^n$. Using phase space methods we
  construct outgoing parametrices and prove Strichartz type estimates
  globally in time. This is done in the context of $C^2$ metrics
  which satisfy a weak asymptotic flatness condition at infinity.  
\end{abstract}


\maketitle

\tableofcontents

\section{Introduction}

Begin with the constant coefficient wave equation in $\R \times \R^n$, $n\ge 2$,
\[
\Box u = (\partial_t^2-\Delta)u = 0, \qquad
u(0)= u_0, \ \ 
u_t(0) = u_1.
\]
On one hand the energy is preserved,
\[
\| \nabla u(t)\|_{L^2} = \|\nabla u(0)\|_{L^2}
\]
where $\nabla $ stands for the space-time gradient of the solution.
On the other hand there is pointwise decay of waves with localized 
initial data. Precisely, if we set $u_0=0$ then
\begin{equation}
\||D_x|^{\frac{1-n}{2}} u(t)\|_{L^\infty} \lesssim 
t^{-\frac{n-1}2} \| u_1\|_{L^1} 
\label{ccdisperse}\end{equation}
for all initial data $u_1$ with a dyadic frequency localization.
As a consequence of this one obtains the Strichartz estimates, which
have the form
\begin{equation}
\||D_x|^{-\rho} \nabla  u\|_{L^pL^q} \leq \|\nabla u_0\|_{L^2} + \|u_1\|_{L^2}.
\label{jg}\end{equation}
This  holds for all pairs $(\rho,p,q)$ satisfying the relations
 $2 \leq p \leq \infty$, $2 \leq q \leq \infty$ and
\begin{equation}
\frac{1}{p}+\frac{n}{q} = \frac{n}{2}-\rho, \ \ \ \ \frac{2}{p}+\frac{n-1}{q}\leq \frac{n-1}{2}  
\label{pq}\end{equation}
with the exception of the forbidden endpoint $(1,2,\infty)$ in
dimension $n=3$.  All $(\rho,p,q)$ satisfying these relations are
called, in the sequel, Strichartz pairs.  If the equality holds in the
second part of (\ref{pq}) then the corresponding pair is called a
sharp Strichartz pair.

In the sequel, we shall not explicitly deal with the case of $q=\infty$,
 allowing us to freely use Littlewood-Paley theory.  When $q=\infty$, 
we only obtain estimates with $L^\infty$ replaced by the appropriate $B_{\infty,2}^{-\rho}$ Besov spaces.  
With additional work, possibly akin to the modifications given to obtain the $\tX$ estimate in \eqref{kxp},
we believe that the proper estimate can be recovered.
In what follows, we shall also concentrate our efforts on the cases that the Strichartz pairs are
sharp.  This can be done without loss of generality since the remaining estimates can then be
recovered using Sobolev embeddings.

A straightforward consequence of (\ref{jg}) is an estimate for
solutions to the inhomogeneous problem
\[
\Box u = f, \ \ \ \ \ \ \ \ \ \ 
u(0) = 0, \ \ \ \ 
u_t(0)= 0,
\]
namely
\begin{equation}
\||D_x|^{-\rho} \nabla u\|_{L^pL^q} \leq \|f\|_{L^1L^2}.
\label{jga}\end{equation}
The simplest case of (\ref{jga}) is the well-known energy estimate
\begin{equation}
\|\nabla u\|_{L^\infty L^2} \leq \|f\|_{L^1L^2}.
\label{jgb}\end{equation}
However, there is a larger family of estimates for solutions to 
the inhomogeneous wave equation where we also vary the 
norms in the right hand side, 
\begin{equation}
\||D_x|^{-\rho} \nabla u\|_{L^p L^q} \leq \||D_x|^{\rho_1} f\|_{L^{p'_1}L^{q'_1}}.
\label{jgc}\end{equation}
This holds for all Strichartz pairs $(\rho,p,q)$, $(\rho_1,p_1,q_1)$.

\def\L{\langle D \rangle}

Estimates of the above type were first proved in the
constant coefficient case in \cite{Brenner}, \cite{Str}.  Further
references can be found in a more recent expository article \cite{MR1151250}.
The endpoint estimate $(p,q)=(2,\frac{2(n-1)}{n-3})$ was only recently
obtained in \cite{MR1646048} ($n \geq 4$).

In this article we are interested in the variable coefficient case of
these estimates, where we replace $\Box$ by a second order 
hyperbolic operator of the form\footnotemark
\[
P(t,x,D) = D_\alpha a^{\alpha \beta}(t,x) D_\beta+ b^\alpha(t,x) D_\alpha + c(t,x).
\]
\footnotetext{Here we employ the summation convention where repeated
  indices are implicitly summed.  Repeated Greek letters
  $\alpha,\beta,\dots$ are summed from $0$ to $n$, where $D_0=D_t$,
  and repeated Latin indices $i,j,\dots$ are summed from $1$ to $n$.}where $D_k = \partial_k / i$.
Here the matrix $a^{\alpha\beta}$ is assumed to have signature $(n,1)$,
and the  time slices are assumed to be space-like, i.e.
$a^{00} < 0$. Thus we consider evolutions of the form
\begin{equation}
Pu = f, \qquad u(0) = u_0, \qquad u_t(0) = u_1.
\label{maineq}\end{equation}

Locally in time this problem is well understood. If the coefficients
are smooth then parametrices are obtained using Fourier integral
operators, and the Strichartz estimates were established in
\cite{MR1168960}.  Operators with $C^{1,1}$ coefficients were first
considered in \cite{MR1644105}, where a wave packet parametrix is
constructed in all dimensions and the Strichartz estimates are proved
in low dimension $n=2,3$.  An alternative parametrix construction,
based on the FBI transform, was later obtained in \cite{nlw}, \cite{cs}
\cite{lp}. There the Strichartz estimates are obtained first for
$C^{1,1}$ coefficients and then for $\nabla^2 a \in L^1 L^\infty$.
Below this regularity threshold for the coefficients the full
Strichartz estimates are lost (see \cite{SS},\cite{MR1909638}), and
one only retains partial results (see \cite{cs},\cite{lp}).

Our goal here is to study the global in time behavior, which is a
considerably more difficult problem. The present article is inspired
by an earlier article of the second author \cite{gS} which deals with
the same issues for the corresponding Schr\"odinger equation.
There are many similarities between the two problems, but also
differences. In what follows we try to discuss both problems
in parallel.

The dynamics for high frequencies are closely related to the Hamilton
flow dynamics, although perhaps less so than in the case of the local
in time problems.  

A first phenomena that one needs to consider is that of refocusing,
which in general precludes the dispersive estimates
\eqref{ccdisperse}  even if we restrict ourselves to coefficients
$a^{\alpha\beta}$ which are sufficiently small, smooth, compactly supported
perturbations of the (Minkowski) identity. This is because even a small
perturbation of the flat metric suffices in order to refocus a group
of Hamilton flow rays originating at the same point and thus produce
caustics.

At the parametrix level this is reflected in the fact that a good
parametrix along a ray which crosses through a bounded region is very
difficult to construct. This is why, following \cite{gS}, we
construct an outgoing parametrix, which only requires 
the analysis of the outgoing Hamilton flow. The price we pay is that 
our parametrix cannot evolve only forward in time; instead it must 
have a forward and a backward component.

 In the case of the Schr\"odinger equation, this is seen on
arbitrarily small time scales due to the infinite speed of
propagation; for the wave equation, on the other hand, a large time
scale is needed.

A second feature is related to the long time behavior of the
bicharacteristics. In the flat case all bicharacteristics are straight
so they escape to infinity both forward and backward in time.
However, in the variable coefficient case, it is possible to have
trapped rays, which are confined to a bounded spatial region.
These correspond to singularities which are largely concentrated
in a bounded region and may destroy not only the dispersive estimates
\eqref{ccdisperse} but also the Strichartz estimates in \eqref{jg}.
On the positive side, the nonexistence of trapped rays is a more stable
phenomena; in particular, it cannot happen for small perturbations. 
Again, this obstruction is seen even on short time
scales for the Schr\"odinger equation, but only on large time scales
for the wave equation.

 The local in time problem for the Schr\"odinger equation has been
 previously considered by other authors.  Stafillani and Tataru~\cite{MR1924470} 
study $C^2$ compactly supported perturbations of the flat metric.
Robbiano and Zuily~\cite{RZ}
 consider smooth asymptotically flat nontrapping metrics in $\R^n$ of
 the short range type and use a parametrix which is a Fourier integral
 operator with complex phase, relying considerably on Sj\"ostrand's
 theory of the FBI transform.  Hassell-Tao-Wunsch \cite{MR2131050}
 have a more direct parametrix construction emulating the
 model of the constant coefficient fundamental solution, which applies
 to smooth asymptotically conic manifolds with short range
 scattering metrics, extended shortly afterward to long range
 scattering metrics.

The dynamics for low frequencies are even more delicate, and for now
there seem to be only two cases where anything at all can be said.
The first is  for sufficiently small perturbations of the flat metric,
which is the case studied in \cite{gS} and here.
The second is for time independent operators, with suitable spectral 
assumptions; for the Schr\"odinger equation this problem is considered
in \cite{MMT} (see also \cite{BT2} and \cite{RT}), while for the wave it will be explored in another forthcoming
paper.

A key part of the global decay estimates are the local energy
estimates, which measure the local averaged decay in the $L^2$
settings.  In the simplest form (see e.g. 
\cite{Alinhac}, \cite{KSS}, \cite{KPV}, \cite{laxphil}, \cite{MetSo}, \cite{morawetz}, 
\cite{smithsogge}, \cite{Strauss}), they are stated
as
\[
\|  \nabla 
u\|_{L^2(\R \times B(0,R))} \lesssim R^\frac12 \|\nabla u(0)\|_{L^2}
\]
when $\Box u = 0$.
Heuristically this is a reflection of the fact that waves move at
speed $O(1)$ and thus spend a time $O(R)$ within a bounded spatial
ball of radius $R$. These are the counterpart of the so called 
local smoothing estimates for the Schr\"odinger equation.  See, e.g.,
\cite{Sjo}, \cite{Veg}, \cite{CoSaut},
\cite{MR1795567}, and \cite{CKS}.
A significant difference is that, in the case of the Schr\"odinger
equation the speed is proportional to the frequency; therefore
one also gains half a derivative in the estimates.

The local energy estimates provide us with a convenient 
space to place the errors in our parametrix and also with a simpler
setup in which to measure the decay of low frequency waves.
In a nutshell, one of our main results asserts that
\[
\text{Local energy estimates } \ \ \Longrightarrow \ \ \ \text{Strichartz estimates}.
\]

The most important part of the article is the outgoing parametrix
construction, for which we are able to adapt the ideas in \cite{gS}.
The parametrix construction in \cite{gS} is based on the use of a time
dependent FBI transform. However it does not use Sj\"ostrand's theory
\cite{Sj1}.  Instead, it takes advantage of the simpler approach
introduced by the second author in \cite{nlw}, \cite{cs}.

For more information about phase space transforms, we refer to
 \cite{MR92k:22017} and \cite{MR93i:35010}.  One of the main starting
 points in the phase space analysis of pde's is Fefferman's article
 \cite{MR85f:35001}. 

Simplified presentations of localized wave packet type parametrix
constructions are now available in \cite{MR2094851},
\cite{phasespace}. These apply to evolutions of the form
\[
(D_t + a^w(t,x,D))u = 0, \qquad u(0) = u_0
\]
on the unit time scale, for symbols $a$ which satisfy
a partial $S^{0}_{00}$ type condition
\[
|\d_x^\alpha \d_\xi^\beta a(t,x,\xi)| \leq c_{\alpha \beta}, \qquad 
|\alpha|+|\beta| \geq 2.
\]
In the finite time analysis in \cite{MR2094851}, \cite{phasespace} the
evolution is turned into a transport equation in the phase space
modulo small errors.  These parametrices are often useful in rescaled
forms. However due to their finite time horizon, they cannot be
directly applied to obtain optimal results for metrics which are not
compactly supported perturbations of the identity.

In the long time analysis in \cite{gS} a time dependent FBI transform
is used instead. A second order term in an asymptotic expansion
becomes nontrivial, and the equation turns into a degenerate parabolic
evolution in the phase space.  Bounds for this evolution are then
obtained using the maximum principle.  Fortunately for us, the main
step in the parametrix construction in \cite{gS} can be applied
directly here for half-waves.  See Theorem \ref{lt1}.

Even though our parametrix is very precise, there are still errors
which need to be controlled and this is done using localized energy
estimates.  We prove such estimates in the case of small perturbations
of the flat metric. For large perturbations nontrapping may fail, and
thus the localized energy estimates may fail. A nontrapping assumption
would help with the localized energy estimates at high frequencies,
but not for the low frequencies.  Here we avoid this problem by using
the localized energy estimates as an assumption for large
perturbations of the flat metric.  In the case of the Schr\"odinger
equation, the local smoothing estimates for large perturbations were considered
 in \cite{MMT}.  See, also, \cite{RT}.  In a follow-up paper we will consider
the same issue in the case of the wave equation.

 Scaling plays an essential role in our analysis. Modulo rescaling and
 Littlewood-Paley theory all our analysis is reduced to waves which
 have fixed frequency of size $O(1)$. Since waves have a propagation
 speed of size $O(1)$,  our study of outgoing waves can be
 largely localized to cones of the form $\{ |x| \approx |t|\}$.
 Certainly the exact flow cannot have a precise localization of this
 type due to the uncertainty principle. To compensate for this we
 introduce an artificial damping term which produces rapid decay of
 waves which do not have the above localization. This allows us to
 restrict our attention to the above cone modulo rapidly decreasing
errors.

 In the present article we consider global in time parametrices and
 Strichartz estimates for $C^{1,1}$ metrics in $\R^n$ which satisfy a
 weak asymptotic flatness assumption.  Due to the global nature of the
 result it is convenient to consider scale invariant assumptions on
 the coefficients.  We denote
\[
A_j = \R \times \{ 2^j \leq |x| \leq 2^{j+1}\}, \qquad A_{<j} =  \R \times \{|x| \leq 2^{j}\}.
\]
Following \cite{gS}, we assume that
\begin{equation}
  \label{coeff}
  \sum_{j \in \Z} \sup_{A_j} |x|^2 |\nabla^2 a(t,x)| + |x| |\nabla a(t,x)|
  + |a(t,x)-M_{1+n}| \leq \e
\end{equation}
where $M_{1+n}$ is the $(n+1)\times (n+1)$ matrix $\text{diag}(-1,1,\dots,1)$
and, for the lower order terms, 
\begin{equation}
  \label{coeffb}
  \sum_{j \in \Z} \sup_{A_j}  |x|^2 |\nabla b(t,x)|
  + |x| |b(t,x)| \leq \e
\end{equation}
\begin{equation}
  \label{coeffc}
  \sup_{\R \times \R^n}   |x|^2 |c(t,x)| \leq \e.
\end{equation}
In some special cases we will need to strengthen the last condition to
\begin{equation}
  \label{coeffcc}
  \sum_{j \in \Z} \sup_{A_j}  |x|^4 |c(t,x)|^2  \leq \e.
\end{equation}
If $\e$ is small enough then \eqref{coeff} precludes the existence of trapped
rays, while for arbitrary $\e$ it restricts the trapped rays to
finitely many dyadic regions.

 Before we state our main results we need to introduce the function
 spaces for the localized energy estimates. We consider a dyadic
 partition of unity in frequency,
\[
1 = \sum_{k=-\infty}^\infty S_k(D_x),
\]
and for each $k \in \Z$ we measure functions of frequency $2^k$ using
the norm
\[
\|u\|_{X_k} = 2^{k/2} \|u\|_{L^2(A_{<{-k}})}+  \sup_{j \geq -k}
\| |x|^{-\frac12}  u\|_{L^2(A_{j})}. 
\]
To measure the regularity of solutions to the wave equation,
we use the  global  norm
\[
\|u\|_{X^s}^2 = \sum_{k=-\infty}^\infty 2^{2sk} \|S_k u\|_{X_k}^2,
\qquad -\frac{n+1}2 < s < \frac{n+1}2.
\]

All Schwartz functions $u \in \mathcal S(\R \times \R^n)$  have finite $X^s$
norm. This allows us to define the space $X^s$ as the completion 
of $\mathcal{S}(\R \times \R^n)$ with respect to the $X^s$ norm. Its structure
is clarified by the next lemma:
\begin{lemma}\cite{gS}
a) ($s=0$) We have
\begin{equation}
\sup_j \| |x|^{-\frac12} u\|_{L^2(A_j)} \lesssim \|u\|_{X^0}. 
\label{hardy0}\end{equation}

b) If $0 < s < \frac{n-1}2$ then the  following  Hardy type inequality holds
for all $u \in \mathcal S(\R \times \R^n)$:
\begin{equation}
\| |x|^{-\frac12-s} u\|_{L^2} \lesssim \|u\|_{X^s}.
\label{hardy1}\end{equation}

c) If $\frac{n-1}{2} \leq s < \frac{n+1}2$ then we have the weaker bound
\begin{equation}
\sum_{j=-\infty}^\infty 2^{-(1+2s)j} \| u - \bar{u}_{A_{<j}}\|_{L^2(A_{<j})}^2
 \lesssim  \|u\|_{X^s}^2
\label{hardy2}\end{equation}
where the time dependent function $ \bar{u}_{A_{<j}}$ stands for the
spatial averages of $u$ in $\{ |x| \leq 2^j\}$.
\label{hd}\end{lemma}
The proof of the lemma is similar to the special case $s=\frac12$
considered in \cite{gS} and is omitted. From the lemma we conclude
that if $s$ is as in case (a,b),  then one can think of $X$ as a space of
distributions.  On the other hand if $s$ is as in case (c),  then $X$ has a BMO type
structure, i.e. $X$ is a space of distributions modulo time dependent
constants.

Controlling the constants is important, particularly when it comes to
localizing parametrices in dyadic regions. This is why we introduce 
also a stronger norm which removes the BMO structure,  namely
\[
\| u \|_{\tX^s}^2 = \| u\|_{X^s}^2 + \| |x|^{-\frac12-s} u\|^2_{L^2},
\qquad 0 < s < \frac{n+1}2.
\]
This coincides with the $X^s$ norm for $0< s < \frac{n-1}2$. To
simplify the exposition we also set $\tX^0 = X^0$.

For the inhomogeneous term in the equation, on the other hand,
we use the dual space $Y^s= (X^{-s})'$ with norm
\[
\|f\|_{Y^s}^2 = \sum_{k=-\infty}^\infty 2^{2sk} \|S_k
f\|_{X_k'}^2, \qquad -\frac{n+1}2 < s < \frac{n+1}{2}.
\]
As $X^s$ is the completion of $\mathcal{S}(\R\times\R^n)$, for $s > -\frac{n+1}2$, the space
$Y^s$ is dense in $S'(\R\times\R^n)$. In addition,
\begin{equation}
\|u\|_{Y^s} \lesssim \||x|^{\frac12-s}  u\|_{L^2}, \qquad \frac{1-n}2
< s < 0
\label{hardy.dual}\end{equation}
and 
\begin{equation}
\|u\|_{Y^0} \lesssim \sum_j \||x|^{\frac12}  u\|_{L^2(A_j)}.
\label{hardy.dual0}\end{equation}

\begin{definition} 
  We say that the operator $P$ satisfies the $\dot{H}^s$ localized energy
  estimates if for each initial data $(u_0,u_1) \in \dot H^{s+1} \times
  \dot H^{s}$ and each inhomogeneous term $f \in L^1 \dot H^{s} + Y^s$,
  there exists a unique solution $u$ to \eqref{maineq}
  with $ \nabla u  \in L^\infty
  \dot H^{s} \cap X^s$ 
which satisfies the bound
\begin{equation}
\|\nabla u\|_{L^\infty \dot H^{s} \cap X^s} \lesssim
\|\nabla u(0)\|_{\dot H^s} + \|f\|_{L^1 \dot H^s + Y^s}.
\end{equation}
\end{definition}

In this context the lower order terms can be often treated as
negligible perturbations:

\begin{lemma}\label{lemma.bc}
a) Let $b$ be as in \eqref{coeffb} and 
\begin{equation}
|s| \leq 1, \qquad |s| <  \frac{n-1}2.
\label{srange}\end{equation}
Then
\begin{equation}
\|b\nabla u\|_{Y^s} \lesssim \e \|\nabla u\|_{X^s}.
\label{bc}\end{equation}

b) Let $n \geq 3$, $c$ be as in \eqref{coeffc} and $-1 < s < 0$. Then 
\begin{equation}
\|c u\|_{Y^s} \lesssim \e \|\nabla u\|_{X^s}.
\label{bcc}\end{equation}

c) Let $n \geq 4$, $c$ be as in \eqref{coeffcc} and $s=-1,0$. Then
\begin{equation}
\|c u\|_{Y^s} \lesssim \e \|\nabla u\|_{X^s}.
\label{bccn4}\end{equation}

d) Let $n = 3$ and $c$ be as in \eqref{coeffcc}. Then
\begin{equation}
\|c u\|_{Y^0} \lesssim \e \|u\|_{\tX^1}.
\label{bccn3}\end{equation}
\end{lemma}

The localized energy estimates hold under the assumption that the 
coefficients $a^{\alpha\beta}$ are a small perturbation of the
Minkowski metric.

\begin{theorem}
   Assume that the coefficients $a^{\alpha\beta}$, $b^\alpha$ satisfy \eqref{coeff},
  \eqref{coeffb} with an $\e$ which is sufficiently small. Assume also
  that $c=0$.  Then the operator $P$ satisfies the $\dot{H}^s$ localized energy
  estimates globally in time for $s$ as in \eqref{srange}.

\label{l2}\end{theorem}

A general coefficient $b$ and a coefficient $c$ can be dealt with
perturbatively but only in dimension $n \geq 3$:

\begin{corollary}
a) Let $n \geq 3$ and $a^{\alpha\beta}$, $b^\alpha$ and $c$ as in \eqref{coeff},
\eqref{coeffb}, \eqref{coeffc} with an $\e$ which is sufficiently small. 
 Then the operator $P$ satisfies the $\dot{H}^s$ localized energy
  estimates globally in time for
\[
-1 < s < 0.
\]

b) Let $n \geq 4$ and $a^{\alpha\beta}$, $b^\alpha$ and $c$ as in \eqref{coeff},
\eqref{coeffb}, \eqref{coeffcc} with an $\e$ which is sufficiently small. 
 Then $P$ satisfies the $\dot{H}^s$ localized energy
  estimates globally in time for  $s=-1,0$.
\label{corc}\end{corollary}

Once we have the local energy estimates, the next step is to construct
an outgoing parametrix which has good time decay and suitable error
bounds in the dual local energy spaces. The parametrix is constructed
at first in the case of a small perturbation of the flat metric.  This
leads to our main scale invariant Strichartz estimate:

\begin{theorem}
  Assume that  $c=0$ and 
the coefficients $a^{\alpha\beta}$, $b^\alpha$ satisfy \eqref{coeff}, \eqref{coeffb} with an
  $\e$ which is sufficiently small.  Let $(\rho_1,p_1,q_1)$ and $(\rho_2,p_2,q_2)$
 be two Strichartz pairs and $s$ as in \eqref{srange}.  Then the solution $u$ to \eqref{maineq}
satisfies
\begin{equation}
\| \nabla u\|_{|D_x|^{\rho_1-s} L^{p_1}L^{q_1} \cap  X^s}
\lesssim \|\nabla u(0)\|_{\dot H^s} + 
\|f\|_{|D_x|^{-\rho_2-s} L^{p'_2}L^{q'_2}+Y^s}.
\label{fse} \end{equation}
A zero order term $c$ can also be added to $P$ subject to the
conditions in  Corollary~\ref{corc}.
\label{tfse}\end{theorem}

If $\e$ is large then any localized energy estimates require an
additional nontrapping condition. Even then the nontrapping can at
most guarantee local in time bounds.  However, we can still prove a
conditional result:

\begin{theorem}
a)  Assume that  $c=0$ and the coefficients $a^{\alpha\beta}$, $b^\alpha$ satisfy
\eqref{coeff}, \eqref{coeffb}.   Then for every
Strichartz pair $(\rho, p,q)$  and $s$ as in \eqref{srange}, we have
\begin{equation}
\||D_x|^{s-\rho} \nabla u\|_{L^p L^q} \lesssim \|\nabla u\|_{ X^s\cap L^\infty\dot{H}^s} + \|Pu\|_{Y^s}.
\label{xpq}\end{equation}
In addition there is a parametrix $K$ for $P$ which satisfies 
\begin{equation}
\|\nabla Kf \|_{ |D_x|^{\rho_1-s} L^{p_1} L^{q_1} \cap X^s} +\|Kf \|_{ \tX^{s+1}}+
 \|(PK -I)f\|_{Y^s}
\lesssim \|f\|_{ |D_x|^{-\rho_2-s} L^{p'_2}L^{q'_2}}
\label{kxp}\end{equation}
for any two Strichartz pairs $(\rho_1,p_1,q_1)$ and
$(\rho_2,p_2,q_2)$.   A  zero order term $c$ can also
be added to $P$ subject to the conditions in Corollary~\ref{corc}.

b) Assume that in addition the operator $P$ satisfies the $\dot{H}^s$ localized
energy estimates. Then the solution $u$ to \eqref{maineq} satisfies
the full Strichartz estimates in \eqref{fse}.
\label{tlargee}\end{theorem}

In applications one might be concerned that the condition
\eqref{coeff} imposes the nontrivial restriction $a(t,0)=M_{1+n}$. This is
true, but it is needed only because we are allowing the derivatives of
the coefficients to be singular at $0$. Otherwise, such a restriction
is unnecessary:

\begin{remark} 
Assume that the the condition \eqref{coeff}  on the
coefficients $a^{\alpha\beta}$ is modified for $|x| < 1$ to 
\[
 \sup_{|x| < 1} (|\nabla^2 a(t,x)| +  |\nabla a(t,x)|
  + |a(t,x)-M_{1+n}|) \leq \e,
\]
and similarly for \eqref{coeffb}, \eqref{coeffc}, and \eqref{coeffcc}.
Assume also that for $k > 0$ the definition of the space $X_k$
is changed to 
\[
\|u\|_{X_k} =  \|u\|_{L^2(A_{<{-0}})}+  \sup_{j \geq 0}
\| |x|^{-\frac12}  u\|_{L^2(A_{j})}.
\]
Then the results in Theorems~\ref{l2},\ref{tfse},\ref{tlargee} remain
valid. Their proofs are  essentially identical with only a few obvious
changes.
\end{remark}

The paper is structured as follows. After introducing some notations
in the next section and making a reduction to the case $a^{00}=-1$ in the third section, 
in Section~\ref{para} we consider the
paradifferential calculus associated to our problem. More precisely,
we show that without any loss we are allowed to mollify the
coefficients $a^{\alpha\beta}$ on a suitable $x$ dependent scale.   This
 allows us to reduce our analysis to problems which are frequency
 localized in dyadic regions.  We also
prove the bound, Lemma \ref{lemma.bc}, for the lower order terms.

Section~\ref{mora} contains the proof of the localized energy
estimates in Theorem~\ref{l2}. The main step of the proof is carried 
out in a frequency localized context and involves a Morawetz type
multiplier technique.

After making a reduction to the half-wave operator in Section 6,
 we state our main result on the existence of
frequency localized outgoing parametrices for half-wave equations,
namely Proposition~\ref{K0} in Section 7.  Using this result we conclude the proof
of Theorems~\ref{tfse},\ref{tlargee}.

The rest of the paper is devoted to the parametrix construction.  
This largely follows \cite{gS}. In
Section~\ref{pdo} we introduce the pseudodifferential operators and
the phase space transforms. An important role is played by the 
 conjugation of pdo's with respect to the phase space
transform, for which we use some results from \cite{MR1944027},
\cite{phasespace}.  In a  first step,  the parametrix is obtained in
Section~\ref{modeleq} in the case of evolutions governed by a
pseudodifferential operator $a^w$ whose symbol satisfies a suitable
smallness condition uniformly in $x$; for this we are fortunately able
to apply directly the result proved in \cite{gS}. This construction is then
transferred in Section~\ref{oureq} to small perturbations of half-waves
via conjugation with respect to the flat half-wave flow.  Finally
to arrive at the desired setup we need to insure that the parametrix
is localized in outgoing propagation cones.  This is done in the last
section by means of choosing a suitable damping term in the equation.

{\em Acknowledgement:} The authors thank the anonymous referee for a
careful reading of the original version of this article and for the
astute suggestions that have greatly enhanced the exposition.

\section{Notations}
\label{notations}

We consider a smooth spatial Littlewood-Paley decomposition
\[
1 = \sum_{j= -\infty}^\infty \chi_j(x) \qquad \supp \chi_j \subset
\{2^{j-1} < |x| < 2^{j+1} \}.
\]
We also set 
\[
\chi_{<j} = \sum_{k<j} \chi_k.
\]

Given $\e$ as in \eqref{coeff}, we can find a sequence  $\e_j \in l^1$
so that 
\begin{equation}
\sup_{A_j}  |x|^2 |\nabla^2 a(t,x)| + |x| |\nabla a(t,x)|
  + |a(t,x)-M_{1+n}|  \leq \e_j
\label{ej}\end{equation}
and 
\[
\sum \e_j \lesssim \e.
\]
Without any restriction in generality,  we can assume that $\e_j$ 
is slowly varying, say
\begin{equation}\label{sv}
|\ln \e_j - \ln \e_{j-1}| \leq 2^{-10}.
\end{equation}
We also choose a function $\e$ in $\R^+$ with the property that
\[
\e_j < \e(s) < 2 \e_j \qquad\text{for} \  2^j < s < 2^{j+1},
\]
and so that 
\[
|\e'(s)|  \leq 2^{-5} s^{-1} \e(s).
\]
This implies that 
\[
\int_0^\infty \frac{\e(s)}{s}\:ds \approx \e.
\]

We also define $\e_k(s)$ so that
$$\e_k(s)\approx \e_j,\quad s\approx 2^j,\quad j\ge -k$$
$$\e_k(s)\approx \e_{-k},\quad s\le 2^{-k}.$$
Note that 
$$\e_k(|x|)\approx \e(2^{-k}+|x|).$$

We consider a frequency Littlewood-Paley decomposition
\[
1 = \sum_{j=-\infty}^\infty S_j(D_x) 
\]
where
\[
\supp s_j \subset \{ 2^{j-1} < |\xi| < 2^{j+1}\}.
\]
We also use the related notations $S_{<k}$, $S_{>k}$, etc.

We say that a function $f$ is localized at frequency $2^k$ if
$\hat{f}$ is supported in $ \{ 2^{k-1} < |\xi| < 2^{k+1}\}$.  An
operator $K$ is localized at frequency $2^k$ if for any $f$  
localized at frequency $2^k$ its image $Kf$ is frequency localized in 
$ \{ 2^{k-10} < |\xi| < 2^{k+10}\}$.

\section{A minor simplification}
\label{minsim}

The aim of this section is to reduce the problem to the case 
when $P$ has the form
\[
P = -D_t^2 +  2D_i a^{i0} D_t + D_i a^{ij} D_j + b^\alpha D_\alpha + c,
\]
and once this is accomplished, $P$ will be taken to be of this form
throughout the sequel.
To arrange that $a^{00} = -1$ we multiply the operator $P$ 
by $- (a^{00})^{-1}$ which satisfies the same bounds as $a^{00}$.
This modifies the other coefficients
\[
a^{\alpha\beta} \to - a^{\alpha\beta} (a^{00})^{-1}, \qquad b^j \to - b^j (a^{00})^{-1}+
D_\alpha((a^{00})^{-1})a^{\alpha j},
\]
\[
b^0\to -b^0(a^{00})^{-1} -(a^{00})^{-1}(D_ta^{00}) +
D_j((a^{00})^{-1})a^{j0},
\qquad c\to -c(a^{00})^{-1},
\]
and it is easy to verify that the assumptions \eqref{coeff},
\eqref{coeffb}, \eqref{coeffc}, and \eqref{coeffcc} are left unchanged. 

To express the second term in the form above we note that
\[
D_t a^{0i} D_i = D_i a^{0i} D_t + (D_t a^{0i}) D_i - (D_i a^{0i}) D_t.
\]
This changes the coefficients $b^\alpha$ but still within the allowed
limits.  Arguing similarly and picking up only lower
order errors within the limits, we may assume that $a^{ij}=a^{ji}$.  
We also note that the coefficient $c$ is not affected by
these transformations.  

To conclude our simplification, we need to verify that our function
spaces are not affected by multiplication by $(a^{00})^{-1}$.

\begin{lemma}
Let $a$ be as in \eqref{coeff} and $s$ as in \eqref{srange}. Then
\begin{equation}
\| a f\|_{Y^s} \lesssim \|f\|_{Y^s}.
\label{afys}\end{equation}
In addition, for all Strichartz pairs $(\rho,p,q)$,  we have 
\begin{equation}
\| a f\|_{ |D_x|^{-\rho-s} L^{p'}L^{q'}} \lesssim \|f\|_{
  |D_x|^{-\rho-s} L^{p'}L^{q'}}.
\label{aflp}\end{equation}
\end{lemma}

\begin{proof}
We write \eqref{afys} in the dual form
\[
|\la a f, u \ra| \lesssim \|f\|_{Y^s} \|u\|_{X^{-s}}
\]
and take a simultaneous Littlewood-Paley decomposition of the three
factors $a$, $f$ and $u$. Nontrivial output is obtained when the two
larger frequencies are comparable. Hence there are three cases to
consider. The trivial one is when the $a$ factor has the low frequency.
For the remaining two cases it suffices to prove the 
off-diagonal decay
\[
|\la S_k a S_k f, S_j u \ra| \lesssim 2^{(s -\delta)(k-j)}  \|S_k
f\|_{X_k'} \|S_j u\|_{X_j}, \qquad j \leq k,
\]
respectively
\[
|\la S_k a S_j f, S_k u \ra| \lesssim 2^{(-s -\delta)(k-j)}  \|S_j
f\|_{X_j'} \|S_k u\|_{X_k}, \qquad j \leq k
\]
for $s$ as in \eqref{srange}. This follows from the definition of the
$X_k'$ and $X_k$ norms combined with the bound on $S_k a$,
\[
|S_k a(x)| \lesssim 2^{-2k}(2^{-2k} +|x|^2)^{-1},
\]
and an uncertainty principle bound for the low frequency factor 
on the dual spatial scale,
\[
\| S_j u\|_{L^2L^\infty(A_{<-j})} \lesssim 2^{\frac{n-1}2 j} \|
S_j u\|_{X_j}, \qquad
 \| S_j f\|_{L^2L^\infty(A_{<-j})} \lesssim 2^{\frac{n+1}2 j} \|S_j f\|_{X_j'}.
\]
The details are straightforward and are left for the reader.

We now prove \eqref{aflp}.  The time variable plays no role and is neglected in the sequel.
We shall use the following variant of a Moser estimate:
\[ \|fg\|_{\dot{W}^{s,p}}\lesssim \|f\|_{L^{q_1}}\|g\|_{\dot{W}^{s,q_2}}+\|g\|_{L^{r_1}}\|f\|_{\dot{W}^{s,r_2}}\]
with $s>0$, $1<p<\infty$, and
\[\frac{1}{p}=\frac{1}{q_1}+\frac{1}{q_2}=\frac{1}{r_1}+\frac{1}{r_2},\quad q_2,r_2\in (1,\infty),
\quad q_1,r_1\in (1,\infty].\]
See, e.g., \cite[\S 2.1, Proposition 1.1, p. 105]{MR1766415}.  We first assume that
$s+\rho\ge 0$ and apply the above estimate to $af$.
This yields
\begin{equation}\label{moser}
\|af\|_{\dot{W}^{\rho+s,q'}} \lesssim \|a\|_\infty \|f\|_{\dot{W}^{\rho+s,q'}}
+ \|a\|_{\dot{W}^{\rho+s,r_1}}\|f\|_{r_2}.
\end{equation}
The first term on the right is trivially bounded by the right side of \eqref{aflp}.  For the second
term on the right, we first pass from the Sobolev space to an appropriate Besov space, and then we
use the following consequence of \eqref{coeff}
\[
| S_l a(t,x)| \lesssim \left\{ \begin{array}{cc}
    2^{-2m -2l}  \e_m &
  2^m <|x| < 2^{m+1}, \ m+l \geq 0 \cr
\e_{-l} & |x| < 2^{-l}.\end{array}\right.\]
Indeed, we see that
\begin{align*}
  \||D_x|^{\rho+s} a\|_{L^{r_1}} &\lesssim \sum_k 2^{(\rho+s) k} \|S_k a\|_{L^{r_1}}\\
&\lesssim \sum_k 2^{(\rho+s)k}\Bigl[ \sum_{m\ge -k} 2^{-2m-2k} \e_m 2^{nm/r_1}
+ 2^{-nk/r_1}\e_{-k}\Bigr].
\end{align*}
Using that $\rho=\frac{n+1}{4}-\frac{n+1}{2q}$ for a sharp Strichartz pair, we have that
$q' < \frac{n}{\rho +s}$.  Thus, we may choose $r_1=\frac{n}{\rho+s}$.  Substituting this in the previous
calculation, noting that $\rho+s<2$ for $\rho$ a part of a
sharp Strichartz pair as above and $s$ as in
\eqref{srange}, and using the summability of $\{\e_m\}$, we have that
\[ \|a\|_{\dot{W}^{\rho+s,r_1}}<\infty,\quad r_1=\frac{n}{\rho+s}.\]
If we now apply Sobolev embeddings to the second factor, we see that the second term in the right
of \eqref{moser} is also bounded by the right side of \eqref{aflp}.  
As we may use a dual argument if $s+\rho<0$, this
completes the proof.
\end{proof}






\section{The paradifferential calculus}
\label{para}

In order to reduce the problem to a frequency localized context and to
simplify the parametrix construction it is convenient to localize the
coefficients in frequency. This is somewhat more complicated than
usual because the frequency localization scale needs to depend on the
spatial scale.

It suffices to work with only the principal part of the operator $P$,
which we denote by
\[
P_a=-D_t^2+2D_i a^{i0}D_t+D_ia^{ij}D_j.
\]
Given a frequency scale $k$ we define the regularized coefficients 
\[
a^{i\beta}_{(k)} =  \delta^{i\beta}
+   \sum_{l < k-4}  (S_{<l} \chi_{<k-2l})   S_l a^{i\beta}. 
\]
Correspondingly we define the mollified operators
\[
P_{(k)} = -D_t^2 + 2 D_i a^{i0}_{(k)} D_t + D_i a^{ij}_{(k)} D_j
\]
which are used on functions of frequency $2^k$. Roughly speaking,
their coefficients are frequency localized in the region
\[
|\xi| \ll 2^k(1+2^k|x|)^{-\frac12}.
\]
 We also introduce a  global mollified operator
\[
\tP =  \sum_{k=-\infty}^\infty P_{(k)} S_k.
\]

  Due to \eqref{ej} and to the fact that the $\e_j$'s are slowly
  varying, it follows that the dyadic parts of the coefficients will
  satisfy the bounds
\begin{equation}
| S_l a^{i\beta}(t,x)| \lesssim \left\{ \begin{array}{cc}
    2^{-2m -2l}  \e_m &
  2^m <|x| < 2^{m+1}, \ m+l \geq 0 \cr
\e_{-l} & |x| < 2^{-l}.
\end{array}\right.
\label{aijl}\end{equation}
This also allows us to obtain bounds on the coefficients 
of $P_{(k)}$, 
\begin{equation}
\begin{split}
&| \d^\alpha (a^{i\beta}_{(k)}(x)-\delta^{i\beta}) | \leq c_\alpha \e_k(|x|)  2^{|\alpha|k}  (1+2^{k}
|x|)^{-|\alpha|}, \qquad |\alpha| \leq 2
\\
&| \d^\alpha a^{i\beta}_{(k)}(x) | \leq c_\alpha \e_k(|x|) 2^{|\alpha|k} 
(1+ 2^k|x|)^{-1-\frac{|\alpha|}2}, \qquad |\alpha| \geq 2.
\end{split}
\label{coeffak}\end{equation}

The main result of this section shows that we can freely replace
$P_a$ by $\tP$ in Theorems~\ref{l2}, \ref{tfse}, \ref{tlargee}(a). It also shows
that at frequency $2^k$ the operators $\tP$ and $P_{(k)}$ are 
interchangeable.

\begin{proposition}
  Assume that the coefficients $a$ satisfy \eqref{coeff}. Then 
\begin{equation}
\|(\tP - P_{(k)})S_l u\|_{X_k'} \lesssim \e \|S_l \nabla u\|_{X_k}, \qquad |l-k|
\leq 2
\label{tamak}\end{equation}
\begin{equation}
\|[P_{(k)},S_k] u\|_{X_k'} \lesssim \e \| \nabla u\|_{X_k}.
\label{comak}\end{equation}
In addition, for $s$ as in \eqref{srange} the following estimate
holds:
\begin{equation}
\| (P_a-\tP) u\|_{Y^s} \lesssim \e \|\nabla u\|_{X^s}.
\label{amta}\end{equation}

\end{proposition}

\begin{proof}

The proof is very similar to the analogous one in \cite{gS}.
We begin with \eqref{amta}, and write 
$$P_a-\tilde{P}=P_{low}+P_{mid}+P_{high}$$
with\footnote{To be explicit with the order of operations,
  the Littlewood-Paley projectors take precedence, followed by
  multiplication and differentiation using right associativity, and finally addition.  Thus, for
  example, the first term in $P_{low}u$ is understood to be
\[\sum_{k=-\infty}^\infty D_i \Bigl[\Bigl(\sum_{l<k-4} (S_{<l}\chi_{\ge k-2l})(S_l a^{ij})\Bigr)\{D_j (S_ku)\}\Bigr].\]
 }
\begin{align*}
  P_{low}&=\sum_{k=-\infty}^\infty D_i \Bigl(\sum_{l<k-4} (S_{<l}\chi_{\ge k-2l})S_l a^{ij}\Bigr)D_j S_k + 
\sum_{k=-\infty}^\infty D_i \Bigl(\sum_{l<k-4} (S_{<l}\chi_{\ge k-2l})S_l a^{i0}\Bigr)D_t S_k
\\
  P_{mid}&=\sum_{k=-\infty}^\infty \sum_{l=k-4}^{k+4}D_i(S_l a^{ij})D_j S_k + 
\sum_{k=-\infty}^\infty \sum_{l=k-4}^{k+4}D_i(S_l a^{i0})D_t S_k
\\
  P_{high}&=\sum_{k=-\infty}^\infty \sum_{l>k+4} D_i(S_l a^{ij})D_j S_k + 
\sum_{k=-\infty}^\infty \sum_{l>k+4} D_i(S_l a^{i0})D_t S_k. 
\end{align*}
Let us examine in detail the second term in each of the expressions above.  The bounds for the 
remaining terms follow from similar arguments.

For $P_{low}$, we notice that the output is at the same frequency $2^k$ as the input.  Since the factor
$D_i$ contributes a factor of $2^k$, it suffices to show
\begin{equation}
  \label{dplow}
\Bigl\|\sum_{l<k-4}(S_{<l}\chi_{\ge k-2l}) S_l a^{i0} v\Bigr\|_{X_k'}\lesssim \e 2^{-k}\|v\|_{X_k}.
\end{equation}

Here, we shall use the bound
\[
|S_{<l}\chi_{\ge k-2l}(x)|\le 
\begin{cases}
  2^{4l-4k},\quad |x|<2^{k-2l-2},\\ 1,\quad |x|\ge 2^{k-2l-2},
\end{cases}
\quad l<k-4.
\]
For $|x|\approx 2^m$, $m\ge -k$, we use this and \eqref{aijl} to see that
\begin{align*}
  \Bigl|\sum_{l<k-4}(S_{<l}\chi_{\ge k-2l}) S_l a^{i0}\Bigr|&\lesssim
\sum_{l=-\infty}^{-m-1} 2^{4l-4k} \epsilon_{-l} + \sum_{l=-m}^{\frac{k-m}{2}-1} \e_m 
2^{-2m-2l}2^{4l-4k} + \sum_{l=\frac{k-m}{2}}^{m-4}\e_m 2^{-2m-2l}\\
&\lesssim 2^{-m-k}\e_m.
\end{align*}
For $|x|<2^{-k}$, the argument is easily modified to give
\[
\Bigl|\sum_{l<k-4} (S_{<l}\chi_{\ge k-2l})S_l a^{i0}\Bigr|\lesssim
\e_{-k},\quad |x|<2^{-k},
\]
which, combined with the previous estimate, yields the desired bound
\eqref{dplow}.

For input frequency $2^k$, $P_{mid}$ permits output frequencies
$2^{h}$ for all $h\le k+4$.  We take  $l=k$ for simplicity of
exposition, and consider separately low and high dimensions.  

In low dimension $n=2,3$  the bound for $P_{mid}$ follows 
from the off-diagonal decay
\[
\|S_h D_i(S_k a^{i0}D_t S_k u)\|_{X_h'}\lesssim \e
2^{\frac{n+1}2(h-k)}
 \|D_t S_k u\|_{X_k},\quad
h\le k+2,
\]
 or more simply,
\begin{equation}
\|S_h(S_k a^{i0} v)\|_{X'_h}\lesssim \e 2^{\frac{n+1}2(h-k)}
2^{-h}\|v\|_{X_k}.
\label{llpl}\end{equation}
Writing
\[
S_h(S_k a^{i0}v) = S_h(\chi_{<-k} S_k a^{i0}v)+\sum_{m\ge -k}
S_h(\chi_m S_k a^{i0}v),
\]
it is sufficient to show that 
\begin{align*}
\|S_h(\chi_m S_k a^{i0} v)\|_{X_h'} &\lesssim \e_m
2^{\frac{n+1}2(h-k)}
2^{-h}\|v\|_{X_k},\quad m\ge -k,\\
\|S_h(\chi_{<-k} S_k a^{i0}v) \|_{X_h'}&\lesssim \e_{-k} 
2^{\frac{n+1}2(h-k)}2^{-h}\|v\|_{X_k}.
\end{align*}
For the former, we apply \eqref{aijl} and see that it suffices to show
\begin{equation}\label{lmk}
\|S_h(\chi_m v)\|_{X_h'}\lesssim 2^{\frac{n-1}2(h-k)} 
2^{k+\frac{3m}{2}}\|v\|_{L^2(|x|\approx 2^m)},\quad m+k\ge 0.
\end{equation}
By interpolating the estimates
\[
\|S_h(\chi_m v)\|_{L^2}\lesssim 
\|\chi_m v\|_{L^2},\quad \|xS_h(\chi_m v)\|_{L^2}\lesssim 2^m\|\chi_m v\|_{L^2},
\quad m+h\ge 0,
\]
we obtain
\begin{equation}\label{fs}
  \|S_h(\chi_m v)\|_{X'_h}\lesssim 2^{\frac{m}{2}}\|\chi_m v\|_{L^2}.
\end{equation}
Recalling that we are in the case when $h < k-2$, this yields
\eqref{lmk} when $m+h\ge 0$. 

 For $m+h<0$ we have improved bounds
\[
\|S_h(\chi_m v)\|_{L^2}\lesssim 2^{\frac{n(m+h)}{2}}\|\chi_m
v\|_{L^2},
\quad
\|x S_h(\chi_m v)\|_{L^2}\lesssim 2^{-h} 2^{\frac{n(m+h)}{2}}\|\chi_m
v\|_2,\quad m+h<0,
\]
which upon interpolation yields 
\[
\|S_h(\chi_m v)\|_{X_h'}\lesssim 2^{\frac{m}{2}}2^{\frac{n-1}{2}(m+h)}\|\chi_m v\|_{L^2},
\]
and implies \eqref{lmk}.  The bound when $\chi_m$ is replaced by
$\chi_{<-k}$ is identical to the $m+h<0$ argument above.

In high dimension $n \geq 4$ the bound \eqref{llpl} is replaced by
\begin{equation}
\|S_h(S_k a^{i0} v)\|_{X'_h}\lesssim \e 2^{2(h-k)}
2^{-h}\|v\|_{X_k}.
\label{llph}\end{equation}
whose proof is similar. The only difference is that now the 
worst case is $m = -h$, whereas in low dimension the 
worst case is $m =-k$ ($n=2$) respectively $-k \leq m \leq -h$
($n=3$).

It remains to consider $P_{high}$ whose output is at frequencies $2^l$
with $l>k$ where $2^k$ is the input frequency. In low dimension 
$n=2,3$ it suffices to show that
\[
\|S_l D_i(S_l a^{i0}D_t S_k u)\|_{X'_l} \lesssim \e 2^{\frac{n-1}2(k-l)} \|D_t
S_ku\|_{X_k},\quad l>k+4
\]
or
\begin{equation}
\|S_l(S_l a^{i0} S_k v)\|_{X'_l}\lesssim \e 2^{\frac{n-1}2(k-l)}
2^{-l}\|S_k v\|_{X_k}, \qquad l > k+4
\label{llplh}\end{equation}
which follows by duality from \eqref{llpl}.

In high dimension $n \geq 4$ instead of \eqref{llplh} we have
\begin{equation}
\|S_l a^{i0} S_k v\|_{X'_l}\lesssim \e 2^{k-l}
2^{-l}\|v\|_{X_k}
\label{llphh}\end{equation}
which is still sufficient except for the endpoint $s=1$. 
At the endpoint we are left with no off-diagonal decay. To compensate 
for that we need a stronger version of \eqref{llphh}, namely
\begin{equation}
\| S_{>k+4} D^2_x a^{i0} S_k v \|_{Y^0}\lesssim \e_{-k} 2^{k} \|v\|_{X_k}.
\label{llphhs}\end{equation}
By \eqref{hardy.dual0} it suffices to show that
\[
\sum_{j=-\infty}^\infty
\| |x|^\frac12 S_{>k+4} D^2_x a^{i0} S_k v \|_{L^2(A_j) }
\lesssim \e_{-k} 2^{k} \|v\|_{X_k}.
\]
But this follows from the bound
\[
|S_{>k+4} D^2_x a^{i0}| \lesssim \e_k(|x|) |x|^{-2},
\]
and in the case that $j<-k$, the following consequence of Bernstein's inequality
\[
\|\chi_{<-k} S_k v\|_{L^2L^\infty}\lesssim 2^{\frac{n-1}{2}k}\|S_k v\|_{X_k}.
\]
The constant $\e_{-k}$ is obtained since the worst case is when $j=-k$,
with exponential decay away from it.

The estimate \eqref{tamak} follows from arguments similar to those
used for $P_{low}$.  To show \eqref{comak}, it would suffice to show
$$\|[S_k, a^{ij}_{(k)}] v \|_{X_k'}\lesssim \e 2^{-k} \|v\|_{X_k}$$
and the equivalent statement with $j=0$, which follows directly from
\eqref{coeffak} with $\alpha=1$.
\end{proof}

In a similar manner we prove the bounds of Lemma \ref{lemma.bc} which show that in
high dimension we can completely dispense with lower order terms.

\begin{proof}[Proof of Lemma \ref{lemma.bc}]
This is again quite similar to the related result from \cite{gS}.

From \eqref{coeffb} we may obtain the following bounds
on the frequency localized pieces of the coefficients
\begin{equation}
|S_k b(x)| \lesssim 2^{-k} \e_k(|x|) (2^{-k} +|x|)^{-2}, \quad |S_{<k} b(x)| 
\lesssim \e_k(|x|) (2^{-k} +|x|)^{-1}.
\label{dyadiccoeff}\end{equation}
To prove \eqref{bc}, we first expand 
\[
b^\alpha\nabla_\alpha u =  \sum_j (S_{<j-4}b^\alpha \nabla_\alpha)S_j u
+ \sum_j\sum_{|k-j|\le 4} (S_k b^\alpha\nabla_\alpha)S_j u
+\sum_j\sum_{k>j+4} (S_k b^\alpha\nabla_\alpha)S_j u.
\]

The easiest case is the first term: the low-high interactions.  Here
the output is at the same frequency range as the input.  Thus, it
would suffice to show
\[
\|(S_{<j-4} b) v\|_{X_j'}\lesssim \e \|v\|_{X_j}.
\]
As the $\e(|x|)$ provides summability, this follows directly from
\eqref{dyadiccoeff}.

The third case is the high-low interactions.  Here, for input at
frequency $2^j$, the output is at frequency $2^k$ with $k>j+4$.  We,
thus, measure the output in $X_k'$.  In low dimension $n=2,3$ it
suffices to show that
\[
\|S_k b S_j v\|_{X_k'}\lesssim  \e 2^{\frac{n-1}2 (j-k)} 
\|S_j v\|_{X_j},\quad j+4<k.
\]
This, however, is just a reformulation of \eqref{llplh} since 
$b$ has exactly the same regularity as $D a$.

In high dimension $n\geq 4$  we have the similar relation
\[
\|S_k b S_j v\|_{X_k'}\lesssim  \e 2^{j-k} 
\|S_j v\|_{X_j},\quad j+4<k
\]
which covers all cases but $s=1$. For $s=1$ we replace this 
with \eqref{llphhs} with $Db$  instead of $D^2 a$.

The remaining case, the high-high interactions, is dual to the
previous case.  This finishes the proof of \eqref{bc}.

Finally, \eqref{bcc}, \eqref{bccn4}, and \eqref{bccn3} follow directly from the embeddings
\eqref{hardy0}, \eqref{hardy1} and their duals \eqref{hardy.dual0},
\eqref{hardy.dual}.
\end{proof}

\section{Localized energy estimates}
\label{mora}

Here we prove Theorem~\ref{l2}.  We can assume that $P$ has the form
in Section~\ref{minsim} with $c=0$. Also due to \eqref{bc} we can take
$b=0$.

The theorem is proved via a positive commutator method. Let
$(\alpha_m)_{m \in \Z}$ be a positive slowly varying sequence with
$\sum \alpha_m = 1$.  Correspondingly we define the space
$X_{k,\alpha}$ with norm
\[
\|u\|_{X_{k,\alpha}}^2 = 2^{k} \| u\|_{L^2(A_{<-k})}^2+ 
 \sum_{j \geq -k} \alpha_j \||x|^{-\frac12}  u\|_{L^2(A_j)}^2
\]
and the dual space
\[
\|u\|_{X'_{k,\alpha}}^2 = 2^{-k} \| u\|_{L^2(A_{<-k})}^2+  \sum_{j \geq -k} \alpha_j^{-1} 
\||x|^{\frac12}  
u\|_{L^2(A_j)}^2.
\]

The key step in the proof of Theorem~\ref{l2} is the following 
frequency localized estimate:

\begin{proposition}
Assume that $\e$ is sufficiently small. Then the bound 
\begin{equation}
\| \nabla u\|_{L^\infty L^2 \cap X_{k,\alpha}} \lesssim \|\nabla u(0)\|_{L^2}
+ \|P_{(k)} u\|_{L^1 L^2 + X'_{k,\alpha}} 
\label{locl2} 
\end{equation}
holds for all functions $u \in L^\infty L^2 \cap X_{k,\alpha}$
localized at frequency $2^k$, uniformly with respect to all slowly
varying sequences $(\alpha_m)$ (as defined in \eqref{sv})  with
\begin{equation}
\sum_{m=-k} ^\infty \alpha_m = 1.
\label{suma}\end{equation}
\label{l2loc}\end{proposition}

\begin{proof}
  By rescaling, the problem reduces to the case when $k = 0$.  We may
  without loss increase the sequence $(\alpha_m)$ so that it remains
  slowly varying and
\begin{equation}
  \label{aam}
  \begin{cases}
    \alpha_0 \approx 1,\\
\sum_{m>0}\alpha_m \approx 1,\\
\e_m \le \epsilon \alpha_{m}.
  \end{cases}
\end{equation}
This is  accomplished by redefining
\[
 \alpha_m := \alpha_m + \frac{\epsilon_m}{\epsilon} + 2^{- 2^{-10} m}.
\]
Since $\epsilon \ll 1$, we can  fix another small parameter $\e \ll \delta \ll 1$ so
  that
\[
\e_m \ll \delta \alpha_{m+\log_2\delta}.
\]

Associating to $(\alpha_m)$ a function $\alpha(s)=\alpha_0(s)$ whose
definition is analogous to that of $\e_0(s)$ in Section
\ref{notations}, we have, from the last property,
\[
\e_0(s)\lesssim \epsilon \alpha(s) \ll \delta\alpha(\delta s).
\]

The proof has three ingredients, the first of which is the classical
energy estimate. Since $\e$ is small it follows that the $\partial_t$
vector field is time-like, and the corresponding energy
\[
E_0(u) = \frac12  \|D_t u\|^2 + \frac{1}{2} \langle a_{(0)}^{ij} D_j u, D_i u \rangle    
\]
is positive definite.  Here and throughout $\langle\,\cdot\,,\,\cdot\,\rangle$ is the
$L^2_x(\R^n)$ inner product.  The time derivative of this energy is
\[
\frac{d}{dt}E_0(u) = \Im\la P_{(0)}u,D_tu\ra + \la (\partial_i a^{i0}_{(0)})D_tu,D_tu\ra
+ \frac{1}{2}\la (\partial_t a^{ij}_{(0)})D_ju,D_iu\ra.
\]

The second component of the proof is a Morawetz-type commutator
estimate.  Let $Q(x,D_x)$ be a spatially self-adjoint operator.  On
time slices we obtain
\begin{multline}\label{mest}
  \frac{d}{dt}\Bigl\{-2\Re\la D_tu, Qu\ra + 2\Re\la a_{(0)}^{j0}D_ju, Qu\ra\Bigr\}
=-2\Im\la P_{(0)}u,Qu\ra + \la i[D_i a^{ij}_{(0)}D_j, Q]u,u\ra
\\+2\Re\la i[a^{i0}_{(0)}D_i,Q]u, D_tu\ra
-2\Re\la (\partial_i a_{(0)}^{i0})D_tu,Qu\ra + 2\Re\la (\partial_t a_{(0)}^{i0})D_i u, Qu\ra.
\end{multline}
The point here is that we seek to choose $Q$ so that the commutator
$[D_i a^{ij}_{(0)} D_j, Q]$ is positive on the characteristic set of
the operator $P_{(0)}$.

Finally, to account for the elliptic region, i.e.  away from the
characteristic set of the operator $P_{(0)}$, we use a Lagrangian
term. Precisely, for a real-valued, time-independent, scalar function
$\psi(x)$, we compute
\[
\begin{split}
  \frac{d}{dt}\Im\la (-D_t+2a^{0j}_{(0)}D_j)u,\psi u\ra =&\ \Re\la
  P_{(0)}u,\psi u\ra
-\Re\la a^{ij}_{(0)}D_ju,\psi D_i u\ra
 +\Im\la a^{ij}_{(0)}D_j u, (\partial_i \psi)
 u\ra\\ &\ +2\Im\la(\partial_t a^{0j}_{(0)}) D_j u,\psi u\ra
-2\Im\la (\partial_j a^{0j}_{(0)})D_t u, \psi u\ra
\\ &\ +\la D_t u, \psi D_t u\ra - 2\Re\la a^{0j}_{(0)} D_ju, \psi D_tu\ra.
\end{split}
\]
We consider two additional small parameters $\delta_0$ and $\delta_1$  so that
\[
\e \ll \delta_1  \ll \delta \ll \delta_0 \ll 1
\]
and define the modified energy 
\[
E(u)= E_0(u)-\delta_0 \Re\la (-D_t+a_{(0)}^{j0}D_j)u,Qu\ra
-\delta_1 \Re\la (-D_t+2a^{0j}_{(0)}D_j)u,i\psi u\ra.
\]
Combining the last three relations we obtain
\begin{equation}
\label{enest}
\begin{split}
  \frac{d}{dt}E(u)&\ + \frac{\delta_0}{2}\la i[D_ia^{ij}_{(0)}D_j,Q]u,u\ra
+ \delta_1 \la D_tu,\psi D_t u\ra \lesssim 
\Im\la P_{(0)}u, (D_t+\delta_0 Q+i\delta_1 \psi)u\ra\\ &\
 + \la |\nabla a_{(0)}|\nabla u,\nabla u\ra
+ \delta_0 |\la i[a_{(0)}^{j0}D_j,Q]u,D_t u\ra| +
 \delta_0 \la |\nabla a_{(0)}||\nabla u|, |Qu|\ra
\\ &\ + \delta_1 \la |a| |\nabla_x u|, \psi|\nabla u|\ra 
+ \delta_1\la |a| |\nabla_x u|,|\nabla \psi||u|\ra
+  \delta_1 \la |\nabla a_{(0)}||\nabla u|,|\psi||u|\ra.
\end{split}
\end{equation}

We choose $Q$ as in \cite{gS}. For convenience its properties are
summarized in the following

\begin{lemma}
There exists an operator $Q$ of the form
\[
Q(x,D_x) =  \delta (D x \phi(\delta |x|) + \phi(\delta |x|) xD)
\]
where $\phi$ has the properties

(i) $\phi(s) \approx (1+s)^{-1}$ for  $s > 0$ and $|\partial^k
\phi(s)| \lesssim (1+s)^{-k-1}$ for $k \leq 4$,

(ii) $ \phi(s) + s \phi'(s) \approx (1+s)^{-1} \alpha(s)$  for  $s > 0$,

(iii) $\phi(|x|)$ is localized at frequency $ \ll 1$,

\noindent and which satisfies the bounds
\[
\| Qu\|_{L^2} \lesssim \|u\|_{L^2}
\]
\[
\| Qu\|_{X_{0,\alpha}} \lesssim \|u\|_{X_{0,\alpha}}
\]
\[
\int_\R \langle i [ D_i a^{ij}_{(0)} D_j,Q] u, u \rangle \:dt \gtrsim \delta \|u\|_{X_{0,\alpha}}^2
\]
for all functions $u$ localized at frequency $1$.
\end{lemma}

The function $\psi(|x|)$ is chosen so that
\[ 
\psi(s)\approx \frac{\alpha(s)}{1+s},\quad |\psi'(s)|\ll \psi(s).
\]
We first note that the above properties of $Q$ and $\psi$ insure that $E$
is positive definite;  specifically
\[ 
E(u)\approx \|\nabla u\|_{L^2}^2
\]
for all functions $u$ at frequency $1$.  Moreover, upon integration in
$t$, we can estimate
\[ 
\int_\R \la D_t u, \psi D_t u\ra \:dt \gtrsim \|D_t u\|^2_{X_{0,\alpha}},
\]
and thus, the integral of the left side of \eqref{enest} is bounded below by
\[
\sup_{t \in \R} E(u)(t) - E(u)(0) + \delta \|\nabla_x u\|^2_{X_{0,\alpha}}
+\delta_1 \|\partial_t u\|^2_{X_{0,\alpha}}.
\]
We now examine the right side of \eqref{enest} after integration in
$t$.  Using \eqref{coeffak}, we have
\[
  \int \la |\nabla a_{(0)}|\nabla u,\nabla u\ra + |\la
  i[a_{(0)}^{j0}D_j,Q]u,D_t u\ra| + \la |\nabla a_{(0)}| |\nabla
  u|,|Qu|+|\psi||u|\ra\,dt \lesssim \epsilon \|\nabla u\|^2_{X_{0,\alpha}}.
\]
Similarly, by our choice of $\psi$, we may find a constant $M>0$ so that
\[
\int C\delta_1 \la |a||\nabla_x u|,\psi |\nabla u|\ra + C\delta_1 \la
|a||\nabla_x u|,|\nabla \psi |u\ra\,dt
\le \frac{\delta_1}{2} \|\partial_t u\|_{X_{0,\alpha}}^2
+ M \delta_1 \|\nabla_x u\|_{X_{0,\alpha}}^2
\]
where $C$ is the implicit constant in \eqref{enest}.

Using these bounds to estimate the right side of \eqref{enest} and
using Cauchy-Schwarz, we obtain
\[
\|\nabla u\|^2_{L^\infty L^2} + \delta \delta_0 \|\nabla_x u\|^2_{X_{0,\alpha}}+ \delta_1 \|\partial_t
u\|^2_{X_{0,\alpha}}\lesssim \|\nabla u(0)\|^2_{L^2} + \delta_1^{-1} \|P_{(0)}u\|^2_{L^1 L^2+X_{0,\alpha}'}
\]
provided, say, $\delta \delta_0 >2M\delta_1$.  This concludes the proof of 
Proposition \ref{l2loc}.
\end{proof}

We conclude now the proof of Theorem~\ref{l2}.  Let $(\beta_m)$ be
another slowly varying sequence with
\[
\sum_{m} \beta_m = 1.
\]
Applying Proposition~\ref{l2loc} with $\alpha_m$ replaced by
$\alpha_m+ \beta_m$ we obtain the bound
\[
\| \nabla u\|_{L^\infty L^2 \cap X_{k,\alpha+\beta}} \lesssim \|\nabla
u(0)\|_{L^2}
+ \|P_{(k)} u\|_{L^1 L^2 + X'_{k,\alpha+\beta}}  
\]
for all $u$ localized at frequency $2^k$.  This implies the weaker
estimate
\[
\|\nabla u\|_{L^\infty L^2 \cap X_{k,\alpha}} \lesssim \|\nabla u(0)\|_{L^2}
+ \| P_{(k)} u\|_{L^1 L^2 + X'_{k,\beta}}.
\]
Since any $l^1$ sequence is dominated by a slowly varying $l^1$
sequence, we can drop the assumption that $\alpha$ and $\beta$ are
slowly varying.  Then we maximize the left hand side with respect to
$\alpha \in l^1$ and minimize the right hand side with respect to
$\beta \in l^1$. This yields
\begin{equation}
\| \nabla u\|_{L^\infty L^2 \cap X_{k}} \lesssim \|\nabla u(0)\|_{L^2}
+ \|P_{(k)} u\|_{L^1 L^2+X'_k}.
\label{xkloc}\end{equation}
For an arbitrary function $u \in X^s$, we apply this bound
to $S_k u$. We have
\[
P_{(k)} S_k u = S_k \tP u + [P_{(k)},S_k] u + S_k(P_{(k)}
-\tP) u.
\]
The last two terms  are frequency localized and can be estimated
by  \eqref{tamak} and \eqref{comak},
\[
\| [P_{(k)},S_k] u + S_k(P_{(k)} -\tP) u\|_{X'_k} \lesssim \e
\sum_{|k-l| \leq 2}\|\nabla S_l u\|_{X_k}.
\]
Then after summation we obtain
\[
\begin{split}
  \| \nabla u\|_{L^\infty \dot{H}^s \cap X^s}^2 &\lesssim \sum_k 2^{2sk} \|\nabla S_k
  u\|^2_{L^\infty L^2 \cap X_k} \\ 
& \lesssim \sum_k\Bigl[ 2^{2sk}\|S_k \nabla u(0)\|_{L^2}^2+ 2^{2sk}\|P_{(k)} S_k
  u\|^2_{L^1 L^2 + X'_k}\Bigr] \\ 
& \lesssim
  \|\nabla u(0)\|_{\dot{H}^s}^2 + \sum_k \Bigl[2^{2sk} \|S_k \tilde{P} u\|^2_{L^1 L^2 + X'_k} 
\\&\qquad\qquad\qquad\qquad\qquad+
 2^{2sk}\| [P_{(k)},S_k] u + S_k(P_{(k)} -\tilde{P}) u\|_{X'_k}^2\Bigr] \\ & \lesssim
 \|\nabla u(0)\|_{\dot{H}^s}^2 + \|\tilde{P}
  u\|^2_{L^1 \dot{H}^s + Y^s} + \e \|\nabla u\|_{X^s}^2 \\ & \lesssim \|\nabla u(0)\|_{\dot{H}^s}^2 +
  \|P_a u\|^2_{L^1 \dot{H}^s +Y^s} + \e \|\nabla u\|_{X^s}^2
\qquad \text{ ( by \eqref{amta})}
\\
&\lesssim \|\nabla u(0)\|_{\dot{H}^s}^2 + \|P u\|^2_{L^1\dot{H}^s + Y^s} +
\e\|\nabla u\|_{X^s}^2
\qquad \text{ ( by Lemma~\ref{lemma.bc})}.
\end{split}
\]
For small $\e$ we can neglect the last right hand side term to obtain
\begin{equation}
  \| \nabla u\|_{L^\infty \dot{H}^s \cap X^s}^2\lesssim  
\|\nabla u(0)\|^2_{\dot{H}^s} + 
\label{dualbd}\|P u\|^2_{L^1 \dot{H}^s + Y^s} 
\end{equation}
which holds in any time interval containing $0$.

Reverting the transformation in Section \ref{minsim}, we see that
without any restriction in generality we can write $P$ in its
self-adjoint divergence form.  Assuming that $b=0$, we may then use a
duality argument to show that for any $f\in L^1\dot{H}^s\cap Y^s$,
there is a $v$ solving
\[ 
Pv=f,\qquad v(0)=v_0,\quad v_t(0)=v_1
\]
with
\[
\|\nabla v\|_{L^\infty \dot{H}^s\cap X^s} \lesssim \|\nabla v(0)\|_{\dot{H}^s}
+ \|f\|_{L^1\dot{H}^s + Y^s}.
\]
Due to \eqref{bc} this extends perturbatively to the case of
nonzero $b$.

By \eqref{dualbd}, this solution is unique, and the proof of Theorem
\ref{l2} is concluded.

\section{The half wave decomposition}\label{halfwave}

In this section we reduce the study of the wave equation 
\eqref{maineq} to the study of two half-wave equations.
We first factor the principal symbol as
\[
-\tau^2+2a^{0j}\tau \xi_j + a^{ij}\xi_i \xi_j = -(\tau+a^+(t,x,\xi))(
\tau+a^-(t,x,\xi))
\] 
where $a^{\pm}$ are $1$-homogeneous in $\xi$ satisfying the
symmetry property
\[
a^-(t,x,\xi) = - a^+(t,x,-\xi)
\]
and are chosen so that $a^+ > a^-$.  The symbols $a^{\pm}$ can be
written down explicitly as
\[a^{\pm}(t,x,\xi) = -a^{0j}(t,x)\xi_j\pm \sqrt{(a^{0j}(t,x)\xi_j)^2 + a^{ij}(t,x)\xi_i\xi_j}.\]
In the sequel, we shall, however, only need the properties listed
above.  This will permit us, in Section \ref{pdo} and beyond, to
free up the $a,b,c$ notation.  There the focus will only be on the
half-wave operators and the symbols $a^{\pm}$.  The notations $a,b,c$
will no longer be reserved for the coefficients of $P$ but will be
used for abstract terms which play the analogous roles.

Mollifying the symbols $a^{\pm}$ 
with respect to $x$ as in Section~\ref{para} we obtain the symbols
$a^\pm_{(k)}(t,x,\xi)$ which we use at frequency $2^k$. We note that 
$a^{\pm}_{(k)}$ are {\bf not} the symbols obtained
from the factorization of the principal symbol of $P_{(k)}$; also 
one cannot define them in this way since algebraic operations
(such as square roots) do not preserve the frequency localization.

We also denote 
\[
l(t,x,\xi) = (a^+(t,x,\xi) - a^{-}(t,x,\xi))^{-1}
\]
and let $l_{(k)}(t,x,\xi)$ be the corresponding regularizations.
We note that $l_{(k)}(t,x,\xi)$ is obtained by regularizing $l(t,x,\xi)$ 
and not by algebraically combining  the symbols
$a^\pm_{(k)}(t,x,\xi)$. 

We are interested in operator properties matching the above 
algebraic properties. We work at frequency $1$, but by rescaling the 
results extend to all dyadic frequencies. 

\begin{proposition}
Define the error operators
\begin{equation}
R^{+} = P_{(0)} + (D_t+A^-_{(0)})(D_t+A^+_{(0)}),
\qquad
R^{-} = P_{(0)} + (D_t+A^+_{(0)})(D_t+A^-_{(0)}).
\label{facp0}\end{equation}
Then for all functions $u$ and $f$ localized at frequency $1$, we have
\begin{equation}
\| R^\pm u \|_{X_0'} \lesssim \| \nabla u\|_{X_0}
\label{rpm}\end{equation}
\begin{equation}
\| \langle x\rangle (L_{(0)} (A^+_{(0)} - A^-_{(0)}) - I) f\|_{L^2_x} \lesssim
\|f\|_{L^2_x}
\label{lapm}
\end{equation}
\begin{equation}
\| [L_{(0)},P_{(0)}] u\|_{X_0'} \lesssim
\|\nabla u\|_{X_{0}}.
\label{comlo}
\end{equation}
\end{proposition}

\begin{proof}
We write $R^+$ in the form
\[
\begin{split}
R^+ =& -i \partial_t A_{(0)}^+(t,x,D)   
-2i (\partial_j a^{j0}_{(0)}) D_t 
- i (\partial_j  a^{jk}_{(0)})D_k
\\ & - (A^- A^+)_{(0)} (t,x,D)+
A_{(0)}^-(t,x,D)  A_{(0)}^+(t,x,D).
\end{split}
\]
The first three terms are easily estimated by \eqref{coeffak}. Consider
the remaining two terms.  For $|\xi|\approx 1$ the symbols
$a^{\pm}(t,x,\xi)$ are smooth and homogeneous in $\xi$. Expanding them
into spherical harmonics we can assume without any restriction in 
generality that both $a^{\pm}$ have the form
\[
a^{\pm}(t,x,\xi) = a^\pm(t,x) h^{\pm}(\xi)
\]
where $a^{\pm}(t,x)$ satisfy bounds similar to the bounds for
$a^{ij}$, namely 
\begin{equation}
|a^\pm(t,x)-a^\pm_\infty|+ |x||\nabla a^\pm(t,x)|+ |x|^2 |\nabla^2 a^\pm(t,x)| \lesssim \e(|x|).
\label{apmbd}\end{equation}
 Then the last two terms in $R^{+}$ have the form
\[
\begin{split}
a^-_{(0)}(t,x)  h^{-}(D)a^+_{(0)}(t,x)  h^+(D) -
(a^-a^+)_{(0)}(t,x)  h^{-}(D) h^+(D) = \qquad \qquad
\\
(a^-_{(0)}(t,x) a^+_{(0)}(t,x) - (a^-a^+)_{(0)}(t,x) ) h^{-}(D) h^+(D)
+ a^-_{(0)}(t,x)  [h^{-}(D),a^+_{(0)}(t,x)]  h^+(D).
\end{split}
\]
The operators $h^\pm$ are bounded in $X_0$ on functions of frequency
$1$. The commutator estimate
\begin{equation}
\label{comapm}
[h^{-}(D),a^-_{(0)}(t,x)]: X_0 \to X_0'
\end{equation}
on frequency $1$ functions follows due to the bound 
\[
|\nabla a^-_{(0)}(t,x)| \lesssim \e_0(|x|) \la x\ra^{-1}.
\]
Hence the estimate \eqref{rpm} is proved if we can show that 

\begin{lemma}
Let $a^{\pm}$ be functions satisfying \eqref{apmbd}. Then
\begin{equation}
|a^-_{(0)}(t,x) a^+_{(0)}(t,x) - (a^-a^+)_{(0)}(t,x) | \lesssim \e_0(|x|)
\la x\ra^{-1}.
\label{amap}\end{equation}
\end{lemma}

\begin{proof}
Without any restriction in generality we can assume that
$a^{\pm}_\infty = 0$. As in the case of the coefficients $a^{ij}$,
the regularized functions have size
\[
| a^\pm_{(0)}(t,x)| \lesssim \e_0(|x|).
\]
We separate the contributions coming from
small $x$ and from large $x$.  The contribution  from small
$x$ decays rapidly at infinity,
\[
| (\chi_{\leq 0} a^\pm)_{(0)}(t,x)| \lesssim \e_0(|x|)\la x\ra^{-N},
\]
and the corresponding part in \eqref{amap} will satisfy a similar
bound. Hence without any restriction in generality
we assume that $a^{\pm}$ are both supported in $A_{\geq 0}$.
This allows us to replace \eqref{apmbd} with a better bound
\[
|a^\pm(t,x)|+ \la x \ra|\nabla a^\pm(t,x)|+ \la x\ra^2 |\nabla^2
a^\pm(t,x)| \lesssim \e_0(|x|).
\]
Using the analogues of \eqref{aijl}, this allows us to estimate
the differences
\[
| a^\pm_{(0)}(t,x) -  a^\pm(t,x)| \lesssim \e_0(|x|) \la x\ra^{-1}
\]
and similarly for their product $a^+ a^-$. The conclusion 
of the lemma follows.
\end{proof}

The proof of \eqref{lapm} is virtually identical, the roles of
$a^{\pm}(t,x,\xi)$ are played by $l(t,x,\xi)$ and
$a^+(t,x,\xi)-a^{-}(t,x,\xi)$.

For \eqref{comlo} we expand $l(t,x,\xi)$ in spherical harmonics
and reduce the problem to the case when 
\[
l(t,x,\xi)=l(t,x) h(\xi)
\]
with $l(t,x)$ satisfying \eqref{apmbd}. Then the proof of \eqref{comlo}
reduces to commutator estimates similar to \eqref{comapm}.
\end{proof}

\section{Parametrices and Strichartz estimates}
\label{tts}

Here we reduce the proof of Theorem~\ref{tfse} to the construction of
a suitable parametrix for $D_t +A^\pm_{(0)}$. Our main result
concerning parametrices is

\begin{proposition}
Assume that $\e$ is sufficiently small. Then there are parametrices
$K_0^\pm$ for $D_t + A_{(0)}^\pm$ which are localized at frequency $1$ 
and have the following properties:

(i) $L^2$ bound:
\begin{equation}
\| K_0^\pm(t,s) \|_{L^2_x \to L^2_x} \lesssim 1,
\label{ees}\end{equation}

(ii) Error estimate: 
\begin{equation}
\begin{split}
 \| (1+|x|)^N (D_t + A_{(0)}^\pm)K_0^\pm(t,s)\|_{L^2_x \to L^2_x}  &\ \lesssim
 (1+|t-s|)^{-N}, \qquad t \neq s,
\\ 
 \| (1+|x|)^N  D_t (D_t + A_{(0)}^\pm)K_0^\pm(t,s)\|_{L^2_x \to L^2_x}  &\ \lesssim
 (1+|t-s|)^{-N}, \qquad t \neq s, 
\end{split}
\label{l2error}\end{equation}

(iii) Jump condition: $K^\pm_0(s+0,s)$ and $K^\pm_0(s-0,s)$ are $S^{0}_{1,0}$
type pseudodifferential operators satisfying
\[
(K^\pm_0(s+0,s) - K^\pm_0(s-0,s)) S_0 = S_0, 
\]

(iv) Outgoing parametrix:
\begin{equation}
\| 1_{\{|x| < 2^{-10} |t-s|\}}  K^\pm_0(t,s)\|_{L^2_x \to L^2_x}  \lesssim (1+|t-s|)^{-N},
\label{outl2}\end{equation}

(v) Pointwise decay:
\begin{equation}
 \| K^\pm_0(t,s) \|_{L^1_x \to L^\infty_x} \lesssim (1 +|t-s|)^{-\frac{n-1}2} . 
\label{ptK0}\end{equation}
\label{K0}\end{proposition}

Here $K^\pm_0$ is defined by
\[
K^\pm_0 f(t) = \int_{-\infty}^\infty K^\pm_0(t,s) f(s) ds.
\]
We leave the proof of this result for later sections, and we show that
it implies Theorems~\ref{tfse},\ref{tlargee}. As an intermediate step
we have the following localized Strichartz estimates for the
parametrix:

\begin{proposition}
  The parametrix $K_0^\pm$ given by Proposition~\ref{K0} has the
  following properties:

(i) (regularity) For any Strichartz pairs $(p_1,q_1)$ respectively 
$(p_2,q_2)$  with $ q_1 \leq q_2$ we have
\begin{equation}
\| K_0^\pm f\|_{L^{p_1}L^{q_1}\cap X_0} \lesssim \|f\|_{L^{p'_2} L^{q'_2}}.
\label{lpbd}\end{equation}

(ii) (error estimate) For any Strichartz pair $(p,q)$ we have
\begin{equation}
\|[(D_t+A^\pm_{(0)})K_0^\pm-1] f\|_{X_0'} 
\lesssim \|f\|_{L^{p'} L^{q'}}.
\label{lperror}\end{equation}
In both \eqref{lpbd} and \eqref{lperror} the function $f$ is assumed
to be localized at frequency $1$.
\label{k0lp}\end{proposition}

The proof is identical to the proof of the similar result in
\cite[Proposition 12]{gS} and is omitted.  The proof of \eqref{lpbd}
follows that of the Strichartz estimates in the constant coefficient
case as it consists of interpolating between \eqref{ees} and
\eqref{ptK0}, using a $TT^*$ argument, and applying the
Hardy-Littlewood-Sobolev inequality.  The error estimate
\eqref{lperror} follows somewhat directly from \eqref{l2error}.

We can use the half-wave parametrices to construct a full wave
parametrix. Precisely we have

\begin{proposition}
  Assume that $\e$ is sufficiently small. Then there is a parametrix
  $K_0$ for $P_{(0)}$ which has the following properties:

(i) (regularity) For any Strichartz pairs $(p_1,q_1)$ respectively 
$(p_2,q_2)$  with $ q_1 \leq q_2$, we have
\begin{equation}
\| \nabla K_0 f\|_{L^{p_1}L^{q_1}\cap X_0} \lesssim 
\|f\|_{L^{p'_2} L^{q'_2}}.
\label{lpbdm}\end{equation}

(ii) (error estimate) For any Strichartz pair $(p,q)$ we have
\begin{equation}
\|(P_{(0)} K_0 -1) f\|_{X_0'} \lesssim \|f\|_{L^{p'} L^{q'}}.
\label{lperrorm}\end{equation}

In all of the above the function $f$ is assumed to be localized at
frequency $1$.
\label{k0lpw}\end{proposition}

\begin{proof}
Our first  approximation for $K_0$ is the operator $K_{00}$ defined by
\[
K_{00}  = L_{(0)} (K_0^+ -K_0^-).
\]
The operator $L_{(0)}$ is bounded in both $ L^{p_1}L^{q_1}$ and $X_0$;
therefore from \eqref{lpbd} we obtain part of \eqref{lpbdm}, namely
\[
\| D_x K_{00} f\|_{L^{p_1}L^{q_1}\cap X_0} \lesssim \|f\|_{L^{p'_2}
  L^{q'_2}}.
\]
We can also bound $D_t K_{00} f$ in $X_0$. We have
\[
D_t K_{00} f  = 
[D_t,L_ {(0)}] (K_0^+ -K_0^-) + 
L_{(0)} D_t  (K_0^+ -K_0^-),
\]
and the first commutator is bounded in $X_0$.
For the second term we
use the $X_0$ bound for $L_ {(0)}$ and write
\[
D_t  (K_0^+ -K_0^-) f = (D_t+A^+_{(0)})  K_0^+ f 
- (D_t+A^-_{(0)})  K_0^- f - A^+_{(0)} K_0^+ f + A^-_{(0)}  K_0^- f.
\]
Now we use \eqref{lperror} and the embedding $X_0' \subset X_0$
for the first two terms and the $X_0$ boundedness of   $A^+_{(0)}$ and
$ A^-_{(0)}$. Summing up we have proved that
\[
\| D_t K_{00} f\|_{ X_0} \lesssim \|f\|_{L^{p'_2}
  L^{q'_2}}.
\]
We still have to estimate $D_t K_{00} f$ in $L^{p_1}L^{q_1}$,
but we postpone this for later.

Next we estimate the error
\[
P_{(0)} K_{00} -1.
\]
The kernel $K_{00}(s,t)$ of $K_{00}$ is smooth in $s,t$ away from the 
diagonal. However, we need to compute its singularity on the diagonal.
Due to the property (iii) in Proposition~\ref{K0} we see that 
the jump of $K_{00}$ on the diagonal vanishes, namely
\[
[K_{00}(t,t)]:=K_{00}(t+0,t) - K_{00}(t-0,t) = 0.
\]
However, the jump of the $t$ derivative of $K_{00}(t,s)$ on the 
diagonal is nontrivial. Precisely, we have
\[
\begin{split} [D_t K_{00}(t,t)] =&\ L_{(0)}[D_t (K_0^+-K_0^-)(t,t)]
  \\
  =&\ L_{(0)}(-A^+_{(0)} [K_0^+(t,t)] +A^-_{(0)} [K_0^-(t,t)])
  \\
  &\ + L_{(0)} ([(D_t + A^+_{(0)})K_0^+(t,t)] -[(D_t +
  A^-_{(0)})K_0^-(t,t)])
 \\
  =&\ L_{(0)}(A^-_{(0)} -A^+_{(0)}) + L_{(0)}
 ([(D_t + A^+_{(0)})K_0^+(t,t)] -[(D_t +  A^-_{(0)})K_0^-(t,t)]).
\end{split}
\]
By \eqref{lapm} the first term is close to the identity, while
the second can be estimated by \eqref{l2error}. For $f$ localized at
frequency $1$ we obtain
\begin{equation}
\| (1+|x|) ([D_t K_{00}(t,t)]-I) f \|_{L^{p'_2} L^2} \lesssim 
\|f\|_{L^{p'_2} L^2}  \lesssim \|f\|_{L^{p'_2} L^{q'_2}}.
\label{koof}\end{equation}

Next we compute
\[
 (P_{(0)} K_{00} -1) f = R_0 f +  ([D_t K_{00}(t,t)]-I) f 
\]
where the first term represents the off-diagonal contribution
and the last term represents the contribution due to the jump of 
$D_t K_{00}(t,s)$ on the diagonal.

We use the factorization \eqref{facp0} for $P_{(0)}$ to compute
the kernel of $R_0$,
\[
\begin{split}
 R_0(t,s)  =&\   L_{(0)}
 P_{(0)}(K^+_0 - K^-_0)(t,s)  - [L_{(0)},P_{(0)}] (K^+_0 - K^-_0)(t,s) 
\\
= &\  - L_{(0)}\left ((D_t + A^{-}_{(0)})
(D_t + A^{+}_{(0)}) K^+_0(t,s) -  (D_t + A^{+}_{(0)})
(D_t + A^{-}_{(0)}) K^-_0(t,s)  \right)
\\ &
+ L_{(0)}\left ( R^+    K^+_0(t,s) -  R^-    K^-_0(t,s)\right) -
  [L_{(0)},P_{(0)}] (K^+_0 - K^-_0)(t,s). 
\end{split}
\]
For the expression on the first line we use the $X_0'$ boundedness
of $L_{(0)}$ and $A^{\pm}_{(0)}$, together with the error estimates in 
\eqref{l2error}. For the $R^{\pm}$ terms we use \eqref{rpm} 
together with \eqref{lpbd} and \eqref{l2error}; the latter is needed
to bound the time derivative $D_t K(t,s)$.
Finally, for the last 
terms we use \eqref{comlo}. Summing up, we obtain
\begin{equation}
\| R_0 f\|_{X_0'} \lesssim \|f\|_{L^{p'}L^{q'}}.
\end{equation}
This is an acceptable error.

The expression
\[
f_1 = -([D_t K_{00}(t,t)]-I) f,
\]
however, is not an acceptable error because it does not yield to a
similar bound of its $X_0'$ norm.  It has better decay at infinity;
therefore we can account for it by setting
\[
K_{0} f = K_{00} f + K_{01} f_1 
\]
where
\[
 K_{01} f_1 (t) = -\frac12 \int_{-\infty}^\infty e^{-|t-s|} f_1(s) ds,
\]
which solves
\[
 \partial_t^2  K_{01} f_1 = f_1 + K_{01}f_1.
\]
Using the bound \eqref{koof} for $f_1$, it is easy to see that
$ K_{01} f_1$ satisfies 
\[
\|\nabla K_{01}f_1\|_{L^{p_1}L^{q_1}\cap X_0}
\lesssim \|f\|_{L^{p_2'}L^{q_2'}}.
\]
On the other hand  the $f_1$ component of the error is replaced by
\[
f_2 = P_{(0)} K_{01} f_1 - f_1 = (2D_i a^{i0}_{(0)} D_t + D_i a^{ij}_{(0)} D_j +I) K_{01} f_1,
\]
which we can estimate by 
\[
\| f_2 \|_{X_0'} \lesssim \|  D_t K_{01} f_1\|_{X_0'} + \|  K_{01}
f_1\|_{X_0'} \lesssim \| \langle x \rangle f_1\|_{L^{p'_2}L^2}
\]
and then apply \eqref{koof}.

The last step of the argument is to prove the $L^{p_1} L^{q_1}$ bound
for $D_t K_{0} f$. We will show that for $u$ at frequency $1$ we have 
\[
\| D_t u\|_{L^{p_1} L^{q_1}} \lesssim \|u\|_{L^{p_1} L^{q_1}} + \|
P_{(0)} u\|_{L^{p'_2} L^{q'_2} + X_0'} 
\]
from which the desired bound follows after an application of
\eqref{lperrorm}.  This would follow from
\[
\| D_t u\|_{L^{p_1} L^{q_1}} \lesssim \|u\|_{L^{p_1} L^{q_1}} + \|
P_{(0)} u\|_{(L^{1} + L^{p_1}) L^{q_1}} 
\]
or equivalently,
\[
\| D_t u\|_{L^{p_1} L^{q_1}} \lesssim \|u\|_{L^{p_1} L^{q_1}} + \| g_1\|_{L^{1}  L^{q_1}} 
+  \| g_2\|_{ L^{p_1} L^{q_1}} , \qquad P_{(0)} u = g_1 + g_2.
\]

The above is easily reduced to the case $g_1 = 0$ by substituting
\[
u := u-v, \qquad v = -\frac12 \int_{-\infty}^\infty e^{-|t-s|} g_1(s)
ds 
\]
since
\[
\| v\|_{(L^1 \cap L^\infty) L^{q_1}} + \| D_t v\|_{(L^1 \cap L^\infty)
  L^{q_1}}  \lesssim   \| g_1\|_{L^{1}  L^{q_1}}.
\]
We are left with proving
\[
\| D_t u\|_{L^{p_1} L^{q_1}} \lesssim \|u\|_{L^{p_1} L^{q_1}} + \|P_{(0)} u\|_{ L^{p_1} L^{q_1}} 
\]
which follows from the interpolation inequality
\[
\| D_t u\|_{L^{p_1} L^{q_1}}^2 \lesssim \|u\|_{L^{p_1} L^{q_1}}\|D_t^2
u\|_{L^{p_1} L^{q_1}}.
\]
\end{proof}

Proposition~\ref{k0lpw} is useful only if $\e$ is small. However,
 a similar result holds even if $\e$ is not small:

\begin{proposition}
 Assume that the coefficients $a^{i\beta}$ satisfy \eqref{coeff}.  
  Then there is a  parametrix $K_0$ for $P_{(0)}$ localized 
at frequency $1$ and which satisfies

(i) (regularity) For any Strichartz pairs $(p_1,q_1)$ respectively 
$(p_2,q_2)$  with $ q_1 \leq q_2$, we have
\begin{equation}
\| \nabla K_0 f\|_{L^{p_1}L^{q_1}\cap X_0} \lesssim \|f\|_{L^{p'_2} L^{q'_2}}.
\label{lpbdl}\end{equation}

(ii) (error estimate) For any Strichartz pair $(p,q)$, we have
\begin{equation}
\|[P_{(0)}K_0-1] f\|_{X_0'} \lesssim \|f\|_{L^{p'} L^{q'}}.
\label{lperrorl}\end{equation}
\label{k0lplarge}\end{proposition}

The proof is identical to the proof of the similar result in
\cite[Proposition 15]{gS}
and is omitted. The idea is that the smallness condition is violated
only on finitely many dyadic spatial regions. 
In \cite{gS} it is argued that a fixed dyadic spatial region can be
partitioned into finitely many cubes on which the smallness holds with
respect to a different coordinate frame. The local parametrices are
then assembled together using a partition of unity.
Alternatively, in a fixed dyadic region the problem of constructing a
parametrix as above can be localized to a similar time scale and then
rescaled into a local problem.

\begin{proof}[Proof of Theorems~\ref{tfse},~\ref{tlargee}]
In what follows we work in a time interval $[T^-,T^+]$, possibly infinite.
By \eqref{tamak} we can replace the operator $P_{(0)}$ by $\tilde{P}$ in
  Propositions~\ref{k0lpw}, ~\ref{k0lplarge}.  Rescaling this result we
  obtain similar parametrices $K_j$ at any dyadic frequency $2^j$.  We
  first assemble these dyadic parametrices and set
\[
K = \sum_{j = -\infty}^\infty K_j S_j.
\]
The properties of $K$ are summarized in the next lemma.

\begin{lemma}
The parametrix $K$ for $P_a$ has the following properties:

(i) (regularity)  For any Strichartz pairs $(\rho_1,p_1,q_1)$ respectively 
$(\rho_2,p_2,q_2)$  with $ q_1 \leq q_2$, we have
\begin{equation}
\|\nabla Kf \|_{|D_x|^{\rho_1-s}L^{p_1}L^{q_1}\cap X^{s}}
\lesssim \|f\|_{|D_x|^{-\rho_2-s}L^{p'_2} L^{q'_2}}.
\label{kf}\end{equation}

(ii) (error estimate) For any Strichartz pair $(\rho, p,q)$, we have
\begin{equation}
\|(P_a K - I) f\|_{Y^{s}} \lesssim \|f\|_{|D|^{-\rho-s}L^{p'} L^{q'}}.
\label{lperror2b}\end{equation}
\label{k}\end{lemma}
 
Part (i) follows directly from the Littlewood-Paley
theory\footnote{As mentioned in the Introduction, the Littlewood-Paley theory with respect to the
  spatial variables cannot be used in dimension $n=2$ for the $L^4
  L^{\infty}$, respectively the $L^{4/3} L^1$ norms.  Here we instead only obtain estimates
  in appropriate $l^2$ Besov spaces.}. Similarly we get
part (ii) but with $\tilde{P}$ instead of $P_a$, since we can write
\[
 \tilde P  K - I = \sum_{j \in \Z} (\tilde P - P_{(j)}) K_j S_j 
+ (P_{(j)} K_j - I) S_j
\]
 However, by
\eqref{amta} we can freely interchange $P_a$ and $\tilde{P}$.  Since
\eqref{bc} allows us to further pass from $P_a$ to $P$ with $c=0$,
this establishes the bounds for the first and last terms in the left side of \eqref{kxp}.

A second step is to use duality to establish an $L^2 \to L^p L^q$
bound.  This establishes \eqref{xpq}, i.e. the first part of
Theorem~\ref{tlargee} (a).

\begin{lemma}
If there  is a parametrix $K$ for $P_a$ as in Lemma~\ref{k} 
and $(\rho, p,q)$ is a Strichartz pair, then 
\begin{equation}
\|\nabla u\|_{|D_x|^{\rho-s}L^p L^q} \lesssim 
\|\nabla u\|_{L^\infty \dot{H}^s \cap X^s} + \|P_a u\|_{Y^s}.
\label{xxp}\end{equation}
 \end{lemma}

 \begin{proof}
   Without any restriction in generality we assume that $T^-$ and
   $T^+$ are finite but prove the bound with constants which are
   independent of $T^+$ and $T^-$.  For $g^\alpha \in |D_x|^{s-\rho} L^{p'}
   L^{q'}$ we use integration by
   parts 
\[
\begin{split}
\int_{T^-}^{T^+} \la \nabla u,g\ra\,dt = &\ \int_{T^-}^{T^+}\la \nabla u, P_a Kg\ra\,dt
-\int_{T^-}^{T^+}\la \nabla u, [P_a K-1]g\ra\,dt \\
= &\ \int_{T^-}^{T^+}\Bigl[ - \la P_a u, \nabla \cdot Kg\ra
-\la \nabla u, [P_a K-1]g\ra 
-2\la \partial_i u, (\nabla a^{i0})\cdot \partial_t Kg\ra
\\&\ -\la \partial_j u, (\nabla a^{ij})\cdot \partial_i Kg\ra
 +2\la \partial_i u, (\partial_t a^{i0})\nabla\cdot Kg\ra
-2\la \partial_t u, (\partial_i a^{i0})\nabla\cdot Kg\ra\Bigr]\,dt
\\&\ + \la \nabla u, \partial_t Kg\ra|_{T^-}^{T^+} + \la \partial_t u,\nabla\cdot Kg\ra|_{T^-}^{T^+}
-2\la a^{i0}\partial_i u, \nabla\cdot Kg\ra|_{T^-}^{T^+}
\\&\ - \la \partial_t u,\partial_t Kg^0\ra|_{T^-}^{T^+} + 2\la a^{i0}\partial_i u,\partial_t Kg^0\ra|_{T^-}^{T^+}
+\la a^{ij}\partial_j u,\partial_i Kg^0\ra|_{T^-}^{T^+}.
\end{split}
\]
Then by \eqref{kf} and \eqref{lperror2b} we obtain
\[
\Bigl|\int_{T^-}^{T^+}\la \nabla u,g\ra\,dt\Bigr|\lesssim \|g\|_{|D|^{s-\rho}L^{p'}L^{q'}}
\Bigl(\|\nabla u\|_{L^\infty\dot{H}^s\cap X^s} + \|P_a u\|_{Y^s}\Bigr).
\]
Here we have also used \eqref{bc} with $b$ replaced by $\nabla a$,
which according to \eqref{coeff} satisfies \eqref{coeffb}.
The conclusion follows.
\end{proof}

Next we prove that the conclusion of Lemma~\ref{k} is also valid 
for $q_1 > q_2$:

\begin{lemma}
The parametrix $K$ in Lemma~\ref{k} also satisfies \eqref{kf} 
when $q_1 > q_2$.
\end{lemma}

\begin{proof}
We repeat the computation in the previous lemma 
with 
\[
u = K f, \qquad g^\alpha \in |D|^{s-\rho_1} L^{p'_1} L^{q'_1}.
\]
All the terms are estimated in the same way except for
\[
\int_{T^-}^{T^+}\la P_a u,\nabla \cdot K g \ra\,dt =
\int_{T^-}^{T^+}\la (P_a K - I)f,  \nabla \cdot K g \ra\,dt 
+ \int_{T^-}^{T^+}\la
f,\nabla \cdot K g\ra\,dt
\]
for which we use \eqref{kf} and \eqref{lperror2b} to estimate
\[
\begin{split}
\Bigl|\int_{T^-}^{T^+}\la P_a u,\nabla \cdot K g \ra\,dt\Bigr| \lesssim &\ \|(P_a K - I)f\|_{Y^s}  
\|\nabla K g \|_{X^{-s}} + 
\| f\|_{|D|^{-s-\rho_2} L^{p'_2} L^{q'_2}} \| \nabla K g\|_{|D|^{s+\rho_2} L^{p_2} L^{q_2}}
\\ 
 \lesssim & \ 
\| f\|_{|D|^{-s-\rho_2} L^{p'_2} L^{q'_2}} 
\| g\|_{|D|^{s-\rho_1} L^{p'_1} L^{q'_1}}.
\end{split}
\]
Then as in the previous lemma we obtain
\[
\begin{split}
\Bigl|\int_{T^-}^{T^+}\la \nabla u,g\ra\,dt\Bigr|&\ \lesssim \|g\|_{|D|^{s-\rho_1}L^{p'_1}L^{q'_1}}
\Bigl(\|\nabla u\|_{L^\infty\dot{H}^s\cap X^s} +
\| f\|_{|D|^{-s-\rho_2} L^{p'_2} L^{q'_2}}  \Bigr)
\\ &\ \lesssim \|g\|_{|D|^{s-\rho_1}L^{p'_1}L^{q'_1}}
\| f\|_{|D|^{-s-\rho_2} L^{p'_2} L^{q'_2}} 
\end{split}
\]
which concludes the proof.
\end{proof}

The bound \eqref{kf} on $Kf$ allows us to estimate $\|Kf\|_{X^{s+1}}$.
However, if $s+1 \geq \frac{n-1}{2}$ then in order to conclude the
proof of \eqref{kxp}, i.e. the remainder of  Theorem~\ref{tlargee} (a), we need to have a bound for the
stronger norm $\|Kf\|_{\tX^{s+1}} $.  This is achieved in the next lemma.

\begin{lemma}
  There is a  parametrix $\tK$ for $P_a$  
 which satisfies:

(i) (regularity)  For any\footnotemark Strichartz pairs $(\rho_1,p_1,q_1)$ respectively 
$(\rho_2,p_2,q_2)$,  we have
\begin{equation}
\|\nabla \tK f \|_{|D_x|^{\rho_1-s}L^{p_1}L^{q_1}\cap X^{s}}
+\| \tK f\|_{ \tX^{s+1}} 
\lesssim \|f\|_{|D_x|^{-\rho_2-s}L^{p'_2} L^{q'_2}}.
\label{tkf}\end{equation}

(ii) (error estimate) For any Strichartz pair $(\rho,p,q)$, we have
\begin{equation}
\|(P_a \tK - I) f\|_{Y^{s}} \lesssim \|f\|_{|D|^{-\rho-s}L^{p'} L^{q'}}.
\label{tlperror2b}\end{equation}
\end{lemma}

\footnotetext{We are again largely ignoring the $L^4L^\infty$, respectively $L^{4/3}L^1$, estimates in $n=2$.}

\begin{proof}
Let $K$ be as in   Lemma~\ref{k}.
  If we think of $Kf$ as the sum of its dyadic pieces which are
  measured in $X_k$, then for $s+1 \geq \frac{n-1}{2}$ we fail to
  obtain a $\tX^{s+1}$ bound for $Kf$ due to the accumulation near the
  origin of the contributions below the uncertainty principle scale
  $\{|x| \lesssim |\xi|^{-1}\}$. To remedy this we attempt to remove
these contributions.

We consider a Schwartz function $\phi$ with  
\[
\phi(0) = 1, \qquad \supp \hat{\phi} \subset \{|\xi| \in [1/2,2]\}
\]
and set $\phi_k(x) = \phi(2^{k} x)$.
In a first approximation we replace the parametrix $K$ with
$(1-T)K$, with $T$ defined by
\[
Tu = \sum_{k = -\infty}^\infty T_k S_k u, \qquad 
T_k u  = u(t,0) \phi_k.
\]
This substitution improves the left hand side of \eqref{tkf}.  We shall show that
\begin{equation}
\|(1- T) u\|_{\tX^{s+1}} \lesssim \|u\|_{X^{s+1}}, \qquad
\frac{n-1}2 \leq s+1 < \frac{n+1}2 
\label{umtu}\end{equation}
\begin{equation}
\|\nabla  T Kf \|_{|D_x|^{\rho_1-s}L^{p_1}L^{q_1}\cap X^{s}}
\lesssim \|  f \|_{|D_x|^{\rho_2-s}L^{p_2'}L^{q_2'}}.
\label{ntu} \end{equation}

For functions $u$ localized at frequency $2^k$, we 
have the fixed time pointwise bound
\[
| u(t,0)| \lesssim 2^{\frac{n-1}2 k} \|u\|_{X_k^0}
\]
which implies that
\[
\| T_k u \|_{X_k} \lesssim  \|u\|_{X_k}.
\]
Here, $X_k^0$ is the spatial part of the $X_k$ norm, i.e., $X_k = L^2 X_k^0$.
Hence we easily obtain
\[
\| Tu \|_{X^{s+1}} \lesssim \| u\|_{X^{s+1}}, \qquad 
\| \nabla Tu \|_{X^s} \lesssim \| \nabla u\|_{X^s}.
\]
The $X^s$ bound of \eqref{umtu} follows immediately, and
in order to obtain the $X^s$ bound of \eqref{ntu}, we then apply \eqref{kf}.
The $L^{p_1}L^{q_1}$ estimate uses a similar argument involving a Bernstein estimate, Littlewood-Paley
estimates, and the bound \eqref{kf}.

To prove the $L^2$ part of the bound
\eqref{umtu}, we take advantage of the fact that $((1-T_k) S_k
u)(t,0)=0$ to obtain the better bound
\[
\sup_j \| {|x|}^{-1-\frac{n}2}   (|x|+2^{-k})^{\frac12+\frac{n}2}     
(1-T_k) S_k u \|_{L^2(A_j)} \lesssim \|S_k u\|_{X_k},
\]
which after summation yields
\[
\| |x|^{-s-\frac32} (1-T) u\|_{L^2} \lesssim \|u\|_{X^{s+1}}, \qquad
\frac{n-1}2 \leq s+1 < \frac{n+1}2 .
\]

Consider now the error estimate for $(1-T)K$. We claim that
\begin{equation}
\| \tP T K f - T f\|_{Y^s} \lesssim \| (\tP K - 1) f \|_{Y^s} +
\|\nabla K f\|_{X^s}.
\end{equation}
It is easily seen that $T_k$ is bounded in $X'_k$; therefore 
$T$ is bounded in $Y^s$. It remains to show that
\[
 \|\tP T u - T \tP u\|_{Y^s} \lesssim \| \nabla u\|_{X^s}
\]
which reduces to
\[
 \|P_{(k)} T_k S_k u - T_k P_{(k)} S_k u\|_{X'_k} \lesssim 
\| \nabla S_k u\|_{X_k}.
\]
After rescaling to $k=0$ this is straightforward. What is important 
is that the second order time derivatives cancel. All the remaining 
terms can be estimated separately.

It remains to consider separately the outstanding error estimate for $TKf$.
This cannot be placed in $Y^s$ because it does not have enough time
integrability. Hence we need to add a correction to the 
parametrix $(1-T)K$ which accounts for this. Our final parametrix
$\tK$ has the form
\[
\tK = (1-T) K + R_T f
\]
where the operator $R_T$  verifies the following properties:
\begin{equation}
\|\nabla  R_T f\|_{|D_x|^{\rho_1-s}L^{p_1}L^{q_1}\cap X^s}   + 
\|R_T f\|_{\tX^{s+1}} \lesssim \|f\|_{|D_x|^{-\rho_2-s} L^{p'_2}L^{q'_2}}
\label{arf}\end{equation}
and 
\begin{equation}
\| (\tP R_T - T)f\|_{Y^s} \lesssim \|f\|_{|D_x|^{-\rho-s}L^{p'}L^{q'}}.
\label{tmdtr}\end{equation}
For $Tf$ we have the representation
\[
T f = \sum_{k\in \Z} \phi_k(x) f_k(t), \qquad f_k(t) = S_k f(t,0).
\]
Then we define 
\[
R_Tf = \sum_{k\in \Z}   \sum_{j\ge k} ( \phi_{j+1}(x) -\phi_j(x)) D_t^{-2}
S^t_{>j}  f_k(t).
\]
Here $S^t_j$ is a Littlewood-Paley decomposition in the time-frequency
variable.  That is, 
\[1=\sum_{j=-\infty}^\infty S_j^t(D_t)\]
with 
\[\supp s_j^t\subset \{2^{j-1}<|\tau|<2^{j+1}\}.\]
The notions of $S^t_{>j}$, $S^t_{\le j}$, etc. are then analogous to those
defined in Section \ref{notations}.  $D_t^{-2}$ denotes the operator
with Fourier multiplier $\tau^{-2}$, where $\tau$ is the frequency
variable dual to $t$.

For $f_k$ we estimate
\[
\| f_k\|_{L^{p'_2}} \lesssim 2^{\frac{nk}{q'_2}} \|S_k
f\|_{L^{p'_2} L^{q'_2}} \lesssim  2^{(-s - \rho_2 +\frac{n}{q'_2})k}  \|S_k
f\|_{|D_x|^{-\rho_2-s} L^{p'_2}L^{q'_2}}.
\]
Since 
\[
- \rho_2 +\frac{n}{q'_2} = \frac{n}2 + \frac{1}{p_2},
\]
we obtain
\begin{equation}
\| f_k\|_{L^{p'_2}} \lesssim  2^{(-s +\frac{n}2 + \frac{1}{p_2})k} 
 \|S_k f\|_{|D_x|^{-\rho_2-s} L^{p'_2}L^{q'_2}}.
\label{fkbd}\end{equation} 
We now proceed to estimate $R_Tf$. By Bernstein's inequality  in time
we have
\[
\begin{split}
\| \nabla ( \phi_{j+1}(x) -\phi_j(x)) D_t^{-2}
S^t_{>j}  f_k(t) \|_{|D_x|^{\rho_1-s}L^{p_1}L^{q_1}} \lesssim
&\ 
2^{-\frac{n}{q_1} j} 2^{-j} 2^{(s-\rho_1)j} 2^{(\frac{1}{p'_2} - \frac{1}{p_1}) j} \| f_k\|_{L^{p'_2}}
\\ = &\ 
2^{-(-s +\frac{n}2 + \frac{1}{p_2})j}   \| f_k\|_{L^{p'_2}}
\\ \lesssim &\ 
2^{(-s +\frac{n}2 + \frac{1}{p_2})(k-j)}   \|S_k
f\|_{|D_x|^{-\rho_2-s} L^{p'_2}L^{q'_2}}.
\end{split}
\]
Since 
\[
s < \frac{n-1}2 \leq \frac{n}2+ \frac{1}{p_2},
\]
it follows that we have off-diagonal decay, while the diagonal
summation is controlled by the Littlewood-Paley theory. This works if
$q_1 \neq \infty$. In the special case $q_1 = \infty$ we also need to
observe that the bump functions $ \phi_{j+1}(x) -\phi_j(x)$
concentrate in different spatial regions; therefore cannot produce
pointwise accumulation.

We continue with the $X^s$ norm:
\[
\begin{split}
\| \nabla ( \phi_{j+1}(x) -\phi_j(x)) D_t^{-2}
S^t_{>j}  f_k(t) \|_{X^s} \lesssim
&\ 
2^{-\frac{n}{2} j}  2^{(s-\frac12)j} 2^{(\frac{1}{p'_2} - \frac{1}{2}) j} \| f_k\|_{L^{p'_2}}
\\ = &\ 
2^{-(-s +\frac{n}2 + \frac{1}{p_2})j}   \| f_k\|_{L^{p'_2}}
\\ \lesssim &\ 
2^{(-s +\frac{n}2 + \frac{1}{p_2})(k-j)}   \|S_k
f\|_{|D_x|^{-\rho_2-s} L^{p'_2}L^{q'_2}},
\end{split}
\]
and the summation works out as before.

The $\tX^{s+1}$ norm is next. Taking advantage of the fact that
$(\phi_{j+1}-\phi_j)(0) = 0$ we compute
\[
\begin{split}
\| |x|^{-s -\frac32} ( \phi_{j+1}(x) -\phi_j(x)) D_t^{-2}
S^t_{>j}  f_k(t) \|_{L^2} \lesssim
&\ 
2^{-\frac{n}{2} j}  2^{(s-\frac12)j} 2^{(\frac{1}{p'_2} - \frac{1}{2}) j} \| f_k\|_{L^{p'_2}}
\\ \lesssim &\ 
2^{(-s +\frac{n}2 + \frac{1}{p_2})(k-j)}   \|S_k
f\|_{|D_x|^{-\rho_2-s} L^{p'_2}L^{q'_2}}
\end{split}
\]
where the restriction $s < \frac{n-1}2 $ insures that the norm on the
left is finite. We still have off-diagonal decay, and for the diagonal
summation we can use 
spatial orthogonality.
This concludes the proof of \eqref{arf}.

For the error estimate \eqref{tmdtr} we split
\[
\tP = - D_t^2 + \tP_1.
\]
The expression $\tP_1 R_T f$ is bounded in the same manner as above.
On the other hand we have
\[
-D_t^2 R_T f - T f = \sum_k   \sum_{j\ge k} ( \phi_{j+1}(x) -\phi_j(x)) 
S^t_{\leq j}  f_k(t),
\]
and for the summand on the right we can use again Bernstein's
inequality with respect to $t$.
 \end{proof}

Now we prove \eqref{fse}. 
If 
\[
P_a u = f+g, \qquad f \in |D_x|^{-\rho_2-s}L^{p'_2} L^{q'_2}, \ g \in Y^s,
\]
then we write
\[
u = Kf + v. \qquad 
\]
We use \eqref{kf} to bound $\nabla Kf$ in $|D_x|^{\rho_1-s}L^{p_1}L^{q_1}\cap X^s$. It remains
to  bound $v$, which solves
\[
P_a v = (1-P_a K)f +g.
\]
In the case of Theorem~\ref{tfse} we use successively \eqref{xxp},
Theorem~\ref{l2}, \eqref{kf}, and \eqref{lperror2b}.  We obtain
\[
\begin{split}
  \|\nabla v\|_{|D_x|^{\rho_1-s}L^{p_1}L^{q_1}} &\lesssim \|\nabla v\|_{L^\infty \dot{H}^s \cap X^s} +
  \|P_a v\|_{Y^s} \\ & \lesssim \|\nabla v(0)\|_{\dot{H}^s} + \|P_a v\|_{Y^s}
  \\ & \lesssim \|\nabla u(0)\|_{\dot{H}^s} + \|\nabla Kf\|_{L^\infty \dot{H}^s} +
  \|(1-P_a K)f\|_{Y^s} + \|g\|_{Y^s} \\ & \lesssim \|\nabla u(0)\|_{\dot{H}^s} +
  \|f\|_{|D_x|^{-\rho_2-s}L^{p'_2} L^{q'_2}} + \|g\|_{Y^s}.
\end{split}
\]
This establishes \eqref{fse} with $P$ replaced by $P_a$.  Using \eqref{bc},
\eqref{fse} then follows.

In the case of Theorem~\ref{tlargee} the argument is similar, but
instead of using Theorem~\ref{l2} we assume that the localized energy
estimates hold.
\end{proof}

\section{Pseudodifferential operators and phase space
transforms} \label{pdo}

Here we tersely introduce the microlocal setup which will be required
in the sequel.  A more detailed exposition can be found in
\cite{phasespace}, \cite{gS}, and the references therein.


Precisely, our initial goal is  to provide a phase-space description of the flow
for a pseudodifferential evolution of the form 
\begin{equation}\label{pdoev}
(D_t + a^w(t,x,D))u=0,\quad u(0)=u_0
\end{equation}
with a real symbol $a$.
We begin by introducing a simpler set-up, which suffices in order to
obtain a short time description of the flow. In terms of symbol classes,
we begin with the standard class $S^{0}_{00}$ of symbols $a$ satisfying
\[
|\partial_x^\alpha \partial_\xi^\beta a(x,\xi)|\le c_{\alpha\beta}, 
\qquad |\alpha|+|\beta|\ge 0.
\]
We also need the following generalizations $S^{(k)}=S^{0,(k)}_{00}$ of the above class,
defined by
\[
|\partial_x^\alpha \partial_\xi^\beta a(x,\xi)|\le c_{\alpha\beta}, 
\qquad |\alpha|+|\beta|\ge k.
\]
For a phase space transform we use the Bargman transform $T$ defined by
\[ 
T u(x,\xi) = c_n \int e^{-\frac{(x-y)^2}{2}}e^{i\xi(x-y)} u(y)\:dy.
\]
This is an isometry from $L^2_x(\R^n)$ to $L^2_{x,\xi}(\R^{2n})$ and thus satisfies
$T^*T=I$. However, $T$ is not an isomorphism; instead, its range consists of functions
which satisfy the Cauchy-Riemann type equation
\begin{equation}\label{cr0}
i\partial_\xi T = (\partial_x - i\xi)T.
\end{equation}
To each pseudodifferential operator $a^w(x,D)$ we associate its phase space
kernel, i.e. the kernel of the conjugated operator $T a^w(x,D) T^*$.
A simple example of the correspondence between the symbol class and the 
phase space kernel is the relation (see \cite[Theorem 1]{phasespace})
\[
a \in S^{(0)} \Leftrightarrow |K((x,\xi),(y,\eta))| \leq c_N (1+|(x,\xi)-(y,\eta)|)^{-N}
\quad \forall \ N \in \N.
\]
This leads to an easy proof of the Calder\'on-Vaillancourt theorem, which asserts that 
the operator $a^w$ is $L^2$ bounded if $a \in S^{(0)}$.

We now turn our attention to the equation \eqref{pdoev} where we
assume $a$ is a real symbol, is in $S^{(2)}$ uniformly in $t\in
[0,1]$, and is continuous in $t\in
[0,1]$.  This suffices in order to guarantee that \eqref{pdoev} is
well-posed in $L^2$. We let $S(t,s)$ denote the evolution operators
corresponding to \eqref{pdoev}; these are all $L^2$ isometries.
We denote by $K(t,s)$ the phase space kernels of $S(t,s)$, i.e. the 
kernels of $T S(t,s) T^*$. It is natural to try to characterize 
the kernels $K(t,s)$ in terms of the Hamilton flow associated to \eqref{pdoev}:
\begin{equation}\label{Hflow}
\begin{cases}
  \dot{x}=a_\xi(t,x,\xi)\\
\dot{\xi}=-a_x(t,x,\xi).
\end{cases}
\end{equation}
The corresponding phase space evolution is denoted by $\chi(t,s)$.
These are canonical transformations in $\R^{2n}$. Furthermore, the
condition $a \in S^{(2)}$ guarantees that $ \chi(t,s)$ are
bilipschitz uniformly with respect to $(t,s) \in [0,1]$.  As it turns
out, the phase space kernel $K(t,s)$ of $S(t,s)$ can indeed be easily
characterized as follows:
\begin{proposition}{\cite[Corollary 7.4]{phasespace}}\label{flatphase}
  Assume that $a$ is a real symbol in $S^{(2)}$ uniformly in $t\in
  [0,1]$.  Then the phase space kernels $K(t,s)$ of $S(t,s)$ satisfy
\[ 
|K(t,x,\xi,s,y,\eta)|\le c_N (1+|(x,\xi)-\chi(t,s)(y,\eta)|)^{-N}.
\]
\end{proposition}

We also have a corresponding Egorov theorem.  For a pdo $q^w(0)$, we define
its conjugate with respect to the flow by
\[
 q^w(t)=S(t,0)q^w(0)S(0,t).
\] 
Then the counterpart of Egorov's theorem in this setting is 
\begin{proposition}{\cite[Proposition 7.6, Proposition 7.7]{phasespace}}\label{conj01}
  Assume that $a$ is a real symbol in $S^{(2)}$ uniformly in $t\in [0,1]$.  
  \begin{enumerate}
    \item[(a.)]  If $q(0)\in S^{(0)}$, then $q(t)\in S^{(0)}$ uniformly in $t$.
    \item[(b.)]  If $q(0)\in S^{(1)}$, then $q(t)\in S^{(1)}$ uniformly in $t$, and
\[ q(t,x,\xi)-q(0)\circ \chi(0,t) \in S^{(0)}.\] 
  \end{enumerate}
\end{proposition}
The counterpart of this result for $q(0)\in S^{(2)}$ is not valid in general.  
However, we can prove it in a special case, which will be useful later.
\begin{proposition}\label{conj02}
  Let $\lambda \geq 1$. Assume that $a(t,x,\xi)=\lambda
  |\xi|$, and let $q(0)\in S^{(2)}$ be an operator which is localized
  at frequency $\lambda$.  Then, $q(t)\in S^{(2)}$ uniformly in $t\in [0,1]$
  and
\[
 q(t,x,\xi)-q(0)\circ\chi(0,t)\in S^{(0)}.
\]
\end{proposition}
 Here,
analogous to the definition in Section 2, we say that an operator $K$
is localized at frequency $\lambda$ if both $\widehat{Kf}$ and
$\widehat{K^*f}$ are supported in
$\{2^{-10}\lambda<|\xi|<2^{10}\lambda\}$ for all $f \in \mathcal{S}(\R^n)$.

We remark that, in the context of the Schr\"odinger equation, a
similar result was proved in \cite{gS} for $a(\xi) = \xi^2$. 

\begin{proof}[Proof of Proposition \ref{conj02}]
  We explicitly compute
\[ 
\chi(t,s)(x,\xi)=(x+\lambda (t-s)\xi |\xi|^{-1},\xi).
\]
Then we want to show that
\[ 
r^w(t,x,D)=e^{-it\lambda |D|}q^w(0,x,D)e^{it\lambda |D|} - q^w(x-\lambda
t D|D|^{-1},D)\in OPS^{(0)} 
\]
uniformly in $t \in [0,1]$.  
Compute
\[
\frac{d}{dt}e^{it\lambda |D|}r^w(t,x,D)e^{-it\lambda|D|} = 
e^{it\lambda |D|}r^w_1(t,x,D)e^{-it\lambda |D|}
 \]
where
\[ 
r_1^w(s,x,D) = -i\lambda [|D|,q^w(x-s\lambda D|D|^{-1},D)] 
- \frac{d}{ds}q^w(x-s \lambda D|D|^{-1},D).
\]
Using the Weyl calculus, as $|\xi| \approx \lambda$ we get
\[ r_1(s,x,\xi)\in S^{(0)}.\]
By Proposition \ref{conj01}, conjugation by $e^{\pm i\lambda t|D|}$ leaves the $S^{(0)}$ class unchanged from which
the conclusion follows.
\end{proof}

From the perspective of the present work, the main disadvantage of
Proposition~\ref{flatphase} is that it can only be used on a fixed
time-scale.  Of course, appropriate versions can be obtained for other
time scales simply by rescaling. For instance, in order to obtain 
results which are valid up to time $s$ we need to replace the Bargman
transform $T$ with its rescaled versions
\[
T_{\frac1s}u(t,x,\xi) = c_n s^{-\frac{n}{4}}\int e^{-\frac{(x-y)^2}{2s}}
e^{i\xi(x-y)} u(t,y)\:dy.
\] 
This is often called the FBI transform. It is still an $L^2$ isometry,
and its range consists of functions satisfying the rescaled
Cauchy-Riemann type equation
\begin{equation}\label{cr}
\frac{i}{s}\partial_\xi T_{\frac{1}{s}} = (\partial_x - i\xi)T_{\frac{1}{s}}.
\end{equation}
Correspondingly, the symbol classes $S^{(k)}$ are replaced by
$S^{(k)}_s$ defined by
\[
a\in S^{(k)}_s:=\{|\partial_x^\alpha \partial_\xi^\beta a(s,x,\xi)|\le 
c_{\alpha\beta}s^{\frac{|\beta|-|\alpha|}{2}},
\quad |\alpha|+|\beta|\ge k\}. 
\]
These are rescaled versions of the $S^{(k)}$ spaces,
and thus, results on $S^{(k)}$ can easily be transferred to
these classes.  In this context, the decay of phase space kernels
would be measured  with a rescaled
distance function
\[ 
d_s((x,\xi),(y,\eta))^2 = s^{-1}|x-y|^2 + s|\xi-\eta|^2.
\]

Still, rescaling does not bring us closer to our goal, which is to
work on an infinite time scale. This difficulty was resolved in
\cite{gS} by using a time dependent scale to study the evolution 
\eqref{pdoev}.

\section{A long time phase space parametrix  }
\label{modeleq}

In this section, following \cite{gS}, we consider global in time evolutions of the form
\begin{equation}
  \label{abevl}
(D_t + a^w(t,x,D)-ib^w(t,x,D)+c^w(t,x,D))u=0,\quad t>0
\end{equation}
with time dependent scales for the symbols $a$, $b$, $c$.
Precisely, we introduce the classes $l^1S^{(k)}$ of symbols in 
 $\R\times T^*\R^n$  whose seminorms
are given by
\[\sum_j 2^{j\bigl(1+\frac{|\alpha|-|\beta|}{2}\bigr)} \|\partial_x^\alpha \partial_\xi^\beta
a(t,x,\xi)\|_{L^\infty(\{t\approx 2^j\})},\quad |\alpha|+|\beta|\ge k. \]
When $k=2$  we also need to better track the second derivatives of the symbols
using the function $\epsilon(t)$ introduced in Section~\ref{notations}.
We denote by $l^1S^{(2)}_\e$ the subset of $l^1S^{(k)}$ whose seminorms are $O(\e)$ when 
$|\alpha|+|\beta|=2$.  This additional condition can be rewritten as
\begin{equation}\label{abbd}
|\partial^\alpha_x \partial^\beta_\xi a(t,x,\xi)|\lesssim 
  \frac{\e(t)}{t}t^{\frac{|\beta|-|\alpha|}{2}},\quad |\alpha|+|\beta|=2.
\end{equation}

We consider the equation \eqref{abevl} 
with a real symbol $a\in l^1S^{(2)}_\e$, which drives the evolution, 
$b\in l^1S^{(1)}$ with $b\ge 0$, which is a damping term, 
and a possibly complex symbol $c\in l^1 S^{(0)}$, which can be regarded as a bounded error. 

We remark that the symbols $a$, $b$, $c$ above are not related to the 
coefficients $a$, $b$, $c$ of $P$, though they play somewhat similar roles.
This slight abuse of notation is harmless since at this stage our arguments
no longer involve the coefficients $a$, $b$, $c$ of $P$; instead all the 
analysis in the parametrix construction is done on the half-wave evolutions
at frequency $1$, using the symbols $a^{\pm}$.

We let $S(t,s)$ now denote the evolution operator corresponding to
\eqref{abevl}. The following result on the $L^2$ evolution was shown in \cite{gS} and follows from fairly standard energy estimate techniques.
\begin{proposition}{\cite[Proposition 28]{gS}}\label{l2tol2}
  Assume that $a \in \SX^{(2)}$ and $b \in \SX^{(1)}$ are real
  symbols with $b \geq 0$, while $c \in \SX^{(0)}$. Then the equation
  \eqref{abevl} is forward well-posed in $L^2(\R^n)$, and the
  corresponding evolution operators satisfy
\[
\|S(t,s)\|_{L^2 \to L^2} \lesssim 1, \qquad 0 < s < t.
\]
\end{proposition}

The evolution \eqref{abevl} is considered in \cite{gS} using a time-dependent
phase space transform. Precisely, at time $t$ one uses the FBI transform
$T_{\frac{1}{t}}$. Thus, the phase space kernels $K(t,s)$ of $S(t,s)$ are  defined to be the kernels of the conjugated operators
\[ 
\tilde{S}(t,s)=T_{\frac{1}{t}}S(t,s)T^*_{\frac{1}{s}}.
\]
A main result in \cite{gS} is to establish precise bounds on the phase
space kernels $K(t,s)$. These  bounds 
are  described in terms of the Hamilton flow dictated by $a$ and the decay dictated by the damping $b$.

The Hamilton flow of $D_t+a^w$ is given by \eqref{Hflow} and as above, we use $\chi(t,s)$ to denote the evolution
operators.  We shall use 
\[ t\to (x_t,\xi_t)\]
to denote the trajectories of the flow.  Using the linearized
equations, one can compute the 
Lipschitz regularity of this flow.  See \cite[Proposition 29]{gS}.  It turns out, however, to be more
convenient to parametrize $\chi(t,s)$ 
using the variables $(x_s,\xi_t)$.  In this context, one obtains the
following regularity.
\begin{proposition}{\cite[Equation (73)]{gS}}
 If $a \in \SX_\e^{(2)}$ with $\e$ small and $s<t$ then
\begin{equation}
  \label{xsxt}
\frac{\partial(x_t,\xi_s)}{\partial(x_s,\xi_t)}=
\begin{pmatrix}
  I_n+\e O(1) & \e O(t)\\
  \e O\bigl(\frac{1}{s}\bigr)&I_n+\e O(1)
\end{pmatrix}.
\end{equation}
\label{chireg}\end{proposition}

In order to describe the decay caused by the damping, we define
\[ 
\psi(t,x_t,\xi_t)=\int_1^t b(s,x_s,\xi_s)\:ds.
\]
We expect $b$ to cause the energy to decay like $e^{-\psi(t,x_t,\xi_t)}$ along the flow.
Using the linearized flow, it can be shown that
\begin{proposition}{\cite[Proposition 30]{gS}}
  If $a \in \SX_\e^{(2)}$ with $\e$ small, $b \in \SX^{(1)}$ and $t >
  s$ then
\begin{equation}
 \frac{\d(\psi(x_t,\xi_t) -\psi(x_s,\xi_s))}{\d(x_s,\xi_t)} = (
 O(s^{-\frac12}), O(t^\frac12)).
\label{psir}\end{equation}
\label{psireg}\end{proposition}

In terms of the above quantities, we can now state the pointwise bound on the kernel of the phase space operator $\tS(t,s)$.  This is one of the principal results of \cite{gS}. 
\begin{theorem}{\cite[Theorem 31]{gS}}
  Let $a \in  \SX_\e^{(2)}$, $b \in  \SX^{(1)}$ be real symbols with
  $b\geq 0$ and $c \in \SX^{(0)}$.  Then for $s < t$ the kernels $K(t,s)$ of the
  operators $\tS(t,s)$ satisfy the bound
\begin{equation}
|K(t,x,\xi_t,s,x_s,\xi)| \lesssim t^{-\frac{n}4} s^\frac{n}4 
\left(1+  (\psi(x_s,\xi_s)-\psi(x_t,\xi_t))^2
+\frac{|x-x_t|^2}t + s|\xi-\xi_s|^2\right)^{-N}.
\label{kbd}\end{equation}
\label{lt1}\end{theorem}

To prove this result one considers the phase space evolution of $
v(t) = T_{\frac{1}t} u(t)$ where $u$ solves \eqref{abevl}.
As it turns out, modulo negligible errors this evolution is 
governed by a degenerate parabolic equation with the following components:

(a) A transport term along the Hamilton flow of $a$

(b) A damping term produced by $b$

(c) A degenerate parabolic term which is due to the change of scale
in the FBI transform.

Pointwise bounds for the kernel of this degenerate diffusion are obtained
in \cite{gS} using the maximum principle.

\section{A perturbation of the half wave equation}
\label{oureq}

The results in Theorem~\ref{lt1} apply for symbols $a$ which satisfy the 
smallness condition $a \in l^1 S_\epsilon^{(2)}$. Instead, the symbols
$a^{\pm}$ are a small perturbation of $\pm |\xi|$, precisely 
\[
a^{\pm} \in \pm |\xi| + l^1 S_\epsilon^{(2)}, \qquad |\xi| \approx 1,
\ |x| \approx t.
\] 
 To remedy this, in this
section we consider the evolution equation
\[
(D_t + |D_x| + a_0^w(t,x,D)-ib_0^w(t,x,D)) u = 0
\]
where $a_0 \in \SX_\e^{(2)}$, $b_0 \in \SX^{(1)}$ are real symbols
with $b_0 \geq 0$.  Since we are interested in this evolution only at frequency
$1$, we will also make the symplifying assumption that
$a_0^w$ is localized at frequency $1$ and that $b_0^w - \tilde{b}_0(t)$ 
is also localized at frequency $1$.  Here $\tilde b_0$ is simply a
function of $t$. 
These assumptions guarantee that if the initial data $u(t_0)$ is localized
at frequencies $\{ |\xi| \in [2^{-10},2^{10}]\}$ then the solution $u$
inherits this localization. The above evolution will serve as the model for our outgoing parametrix.

We denote by $S_0(t,s)$ the $L^2$ evolution generated by
the above equation. Due to the above frequency localization of 
$a_0^w$ and $b_0^w - \tilde{b}$ we have
\[
 S_0(t,s) S_{-10<\cd<10} = S_{-10<\cd<10} S_0(t,s) =  S_{-10<\cd<10} S_0(t,s) S_{-10<\cd<10}.
\]
We denote by $\tS_0(t,s)$ its (frequency localized) phase space image
\[
\tS_0(t,s) = T_{\frac1t} S_0 (t,s) S_{-10<\cd<10} T^*_\frac{1}s.
\]
We want to obtain bounds on the kernel of $\tS_0(t,s)$ which are
similar to the ones in Theorem~\ref{lt1}. As a preliminary step we
need to study the regularity of the associated Hamilton flow
which we denote by $\chi_0(t,s)$. This can be done directly, but 
for our purposes it is more convenient to reduce it to the case
considered in the previous section. 

At each time $t$ we consider the symplectic map $\mu$ defined by
\[
\mu_t(x,\xi) = (x+t\xi|\xi|^{-1},\xi).
\]
which corresponds to the Hamilton flow for the $D_t+|D|$ evolution.
This extends to a space-time symplectic map
\[
\mu(t,\tau,x,\xi) = (t,\tau-|\xi|,x+t\xi|\xi|^{-1},\xi).
\]

If $p_0$ is the symbol 
\[
p_0(t,\tau,x,\xi) = \tau+|\xi|+ a_0(t,x,\xi),
\]
then its image through $\mu$ is
\[
p_0(\mu(t,\tau,x,\xi)) = \tau+ a(t,\tau,x,\xi), \qquad a(t,x,\xi) =
a_0(t,x+t\xi|\xi|^{-1},\xi).
\]
Hence the conjugate of the Hamilton flow $\chi_0(t,s)$ for $\tau+|\xi|
+ a_0$ with respect to $\mu_t$ is the Hamilton flow $\chi(t,s)$ for
$\tau+a(t,x,\xi)$,
\[
\chi_0(t,s) = \mu_t \circ \chi(t,s) \circ \mu_s^{-1}.
\]

We note that $a \in l^1S^{(2)}_\e$ iff $a_0 \in l^1S^{(2)}_\e$. Hence 
from \eqref{xsxt} we obtain its counterpart for the $\chi_0$ flow,

 \begin{proposition}
   If $a_0 \in \SX_\e^{(2)}$ with $\e$ sufficiently small and $t > s$
   then the Hamilton flow $\chi_0(t,s)$ has the Lipschitz regularity
\begin{equation}
 \frac{\d( x_t,\xi_s) }{\d (x_s,\xi_t)}  = 
\left(\begin{array}{cc} I_n+\e O(1) 
&  2(t-s)|\xi|^{-3}  (|\xi|^2 I_n -\xi\otimes \xi)  + \e O(t) 
\cr\e O(\frac{1}s) &
 I_n + \e O(1)  
 \end{array}   \right).
\label{xsxt1}\end{equation}
\label{flowreg1}\end{proposition}

We proceed in a similar manner with $b_0$ and set
\[
b(t,x,\xi) =b_0(t,x + t\xi|\xi|^{-1},\xi).
\] 
Then the integral $\psi_0$ of $b_0$ along the $\chi_0$ flow is the
$\mu$ conjugate of the integral $\psi$ of $b$ along the $\chi$ flow.
Hence we also trivially obtain the analog of Proposition~\ref{psireg},
namely

\begin{proposition}
 If $a_0 \in  \SX_\e^{(2)}$ with $\e$ sufficiently small and $b_0  \in
 \SX^{(1)}$ then for $t > s$ we have
\begin{equation}
 \frac{\d(\psi_0(x_t,\xi_t) -\psi_0(x_s,\xi_s))}{\d(x_s,\xi_t)} = (
 O(s^{-\frac12}),  O(t^\frac12)).
\label{psir1}\end{equation}
\label{psireg1}\end{proposition}

Now we can state our main result:

\begin{theorem}
  Let $a_0 \in \SX_\e^{(2)}$, $b_0 \in \SX^{(1)}$ be real symbols with
  $b_0 \geq 0$ with $\e$ sufficiently small, so that $a_0$ and $b_0 -\tilde b_0(t)$ are localized at frequency $1$. Then for $s < t$ the
  kernel $K_0$ of the operator $\tS_0(t,s)$ satisfies the bound
\begin{multline}
|K_0(t,x,\xi_t,s,x_s,\xi)| \lesssim  t^{-\frac{n}4} s^\frac{n}4 
\left(1+  (\psi_0(x_s,\xi_s)-\psi_0(x_t,\xi_t))^2
+\frac{|x-x_t|^2}t + s|\xi-\xi_s|^2\right)^{-N}\\
\times (1+ t\ d(|\xi_s|,[2^{-10},2^{10}]))^{-N}.
\label{kbd1}
\end{multline}
\label{tphase1}\end{theorem}

\begin{proof}
We use  Theorem~\ref{lt1} via a conjugation with respect to 
the flat half-wave flow, which corresponds to the canonical
transformations $\mu_t$. Denote
\[
S(t,s) = e^{ i t |D|} S_0(t,s) e^{-is |D|}.
\]
Then we compute
\[
\frac{d}{dt} S(t,s) = -ie^{ i t |D|} (-a^w_0(t,x,D) +ib^w_0(t,x,D)) e^{-i
  t |D|} S(t,s).
\]
Hence $S(t,s)$ is the evolution associated to the pseudodifferential
operator
\[
 e^{ i t |D|} (-a^w_0(t,x,D)+ib^w_0(t,x,D))  e^{-i t |D|}. 
\]
Using rescaled versions of Propositions~\ref{conj01},\ref{conj02},
this operator can be expressed in the form
\[
a^w(t,x,D) - ib^w(t,x,D) + c^w(t,x,D)
\]
where the remainder term satisfies $c \in \SX^{(0)}$. Hence the phase
space kernel of $S(t,s)$ satisfies the bounds given by
Theorem~\ref{lt1}.

Returning to the original equation, for the phase space evolution
$\tS_0(t,s)$ we can write
\[
\begin{split}
\tS_0(t,s) &= T_{\frac{1}t} e^{-it|D|} S_{-10<\cd<10} S(t,s) S_{-10<\cd<10}  e^{is|D|}  T^*_{\frac{1}s}
\\ &= T_{\frac{1}t} e^{-it|D|}  S_{-10<\cd<10}  T_{\frac{1}t}^* T_{\frac{1}t}S(t,s) 
 T_{\frac{1}s}^* T_{\frac{1}s}S_{-10<\cd<10} e^{is|D|}  T_{\frac{1}s}^*
\\ &= (T_{\frac{1}t} e^{-it|D|}  S_{-10<\cd<10}  T_{\frac{1}t}^*) \tS(t,s) 
 ( T_{\frac{1}s}S_{-10<\cd<10} e^{is|D|}  T_{\frac{1}s}^*).
\end{split}
\]
By a rescaled version of Proposition~\ref{flatphase} the kernel of the
first factor $T_{\frac{1}t} e^{-it|D|}  S_{-10<\cd<10} T_{\frac{1}t}^*$ is rapidly
decreasing on the $t^{\frac12} \times t^{-\frac12}$ scale away from
the graph of $\mu_t$ as well as away from the support of the symbol
$S_{-10<\cd<10}$, 
while the kernel of the last factor
$T_{\frac{1}s}e^{is|D|} S_{-10<\cd<10} T_{\frac{1}s}^*$ is rapidly decreasing on the
$s^{\frac12} \times s^{-\frac12}$ scale away from the graph of
$\mu_s^{-1}$ as well as away from the support of the symbol $S_{-10<\cd<10}$.  Hence the composition simply replaces the Hamilton flow
associated to $a$ by the Hamilton flow associated to $a_0$ and the
function $\psi$ with $\psi_0$ in the kernel bounds. Thus \eqref{kbd}
implies \eqref{kbd1}, and the proof is concluded.
\end{proof}


\section{The parametrix construction}

We end with a proof of Proposition \ref{K0}.  That is, we construct a parametrix $K_0^+$ for
$D_t+A^+_{(0)}$, and easy modifications yield also a parametrix for $D_t+A^-_{(0)}$.  In the sequel, we shall
drop the $+$ signs and denote these by $K_0$ and $D_t+A_{(0)}$ respectively.  The $\pm$ signs will be reserved
to distinguish waves which are outgoing forward, respectively backward, in time.

We partition the annulus $|\xi|\approx 1$ in phase space
\[ s_{-1}(\xi)+s_0(\xi)+s_1(\xi)=\sum_{\pm}\sum_{j\ge 0} p_j^{\pm}(x,\xi),\]
with 
\[\supp p_j^{\pm}\subset \{2^{j-1}<|x|<2^{j+1},\: \pm x\xi\ge
-2^{-5}|x|\},\quad j\ge 1,\]
\[\supp p_0^{\pm}\subset \{|x|<2,\: \pm x\xi\ge
-2^{-5}|x|\}.\]
At the expense of Schwartz tails which play no role in the sequel, we may replace $p_j^{\pm}$ by
$S_{<-10}(D_x)p_j^{\pm}$.  As such, we shall do so without changing the notation.  This allows us
to assume that the operators $P_j^{\pm}$ are frequency localized to frequency $1$.

In the proposition which follows, we construct evolution operators $S_{j}^{\pm}(t,s)$ as the evolutions associated to a certain
damped half-wave equation.  We then form $K_0$ by setting
\[
K_0(t,s)=
\begin{cases}
  \sum_{j=1}^\infty S_j^-(t,s)(P_j^-)^w(x,D),\quad t<s\\
  \sum_{j=1}^\infty S_j^+(t,s)(P_j^+)^w(s,D),\quad t>s.
\end{cases}
\]
The properties of $K_0$ listed in Proposition \ref{K0} follow easily, after summing, from the
given properties of $S_j^{\pm}$.

\begin{proposition}\label{Kj}
  Assume that $\e$ is sufficiently small.  Then for each $s\in \R$, there is an outgoing parametrix
$S_j^+$ for $D_t+A_{(0)}$ in $\{t>s\}$ which is localized at frequency $1$ and satisfies the following:
  \begin{enumerate}
    \item[(i)] $L^2$ bound:
       \[ \|S_j^+ (t,s)\|_{L^2\to L^2}\lesssim 1 \]
    \item[(ii)] Error estimate:
      \begin{equation}
	\label{l2errorj} 
	\begin{split}
	\|x^\alpha (D_t+A_{(0)})S_j^+(t,s)P_j^+\|_{L^2\to L^2}&\lesssim (2^j + |t-s|)^{-N}\\
        \|x^\alpha D_t(D_t+A_{(0)})S_j^+(t,s)P_j^+\|_{L^2\to L^2}&\lesssim (2^j+|t-s|)^{-N}
	\end{split}
      \end{equation}
    \item[(iii)] Initial data:
        \[ S_j^+(s+0,s)=I\]
    \item[(iv)] Outgoing parametrix:
      \begin{equation}
	\label{outl2j}
	   \|\1_{\{|x|<2^{-10}(|t-s|+2^j)\}} S_j^+(t,s)P_j^+\|_{L^2\to L^2}\lesssim (|t-s|+2^j)^{-N}
      \end{equation}
    \item[(v)] Finite speed:
      \begin{equation}
	\label{finitel2j}
	   \|x^\alpha \1_{\{|x|>2^{10}(|t-s|+2^j)\}} S^+_j(t,s)P_j^+\|_{L^2\to L^2}\lesssim (|t-s|+2^j)^{-N}
      \end{equation}
   \item[(vi)] Frequency localization:
\begin{equation}
	\label{freqloc}
\|(1-P_{[-4,4]}) S^+_j(t,s)P_j^+\|_{L^2\to L^2}\lesssim (|t-s|+2^j)^{-N}
\end{equation}
    \item[(vii)] Pointwise decay:
      \begin{equation}\label{ptwisedecay}
	  \|S_j^+(t,s)P_j^+\|_{L^1\to L^\infty}\lesssim (1+|t-s|)^{-\frac{n-1}{2}}.
      \end{equation}
  \end{enumerate}
With obvious modifications, the same hold for $S_j^-$.
\end{proposition}

By translation invariance, without any loss of generality we  may assume that $s=2^j$. We first reduce the problem to the study of an evolution of a perturbed half-wave equation as in Section \ref{oureq}. Heuristically we 
observe that in the support of the symbol of $P_j^+$ we have 
\[
 |\xi| \in [2^{-2}, 2^2], \qquad |x| \approx s, \qquad x \cdot \xi \geq
-\frac15 |x| |\xi|.
\]
An easy computation shows that  along the forward Hamilton flow starting here we have
\[
 |\xi| \in [2^{-3}, 2^3], \qquad |x| \approx t.
\]
But in this region we have 
\[
a_{(0)}(t,x,\xi)-|\xi|\in l^1S^{(2)}_\e 
\]
which follows from the analog of \eqref{coeffak} which
holds for $a_{(0)}$.  Thanks to \eqref{outl2j}, \eqref{finitel2j}
and \eqref{freqloc}, we can freely modify the symbol of $a_{(0)}$ 
in the regions $\{ |x| \ll t\}$ and $\{|\xi| \not\in [2^{-5},2^5]\}$
at the expense of  producing a negligible error in \eqref{l2errorj}.  
 It thus suffices to study the evolution governed by a symbol
\[
|\xi|+a_0(t,x,\xi),\quad a_0\in l^1S^{(2)}_\e 
\]
so that $a_0$ vanishes if $\{|\xi| \not\in [2^{-6},2^6]\}$ and $a_0^w$ is localized at frequency $1$.

In essence, $a_0=a_{(0)}-|\xi|$, and thus by \eqref{coeffak}, we may assume the better regularity
\begin{equation}
  \label{extraa}
  \begin{split}
    |\partial_x^\alpha \partial_\xi^\beta a_0(t,x,\xi)|&\lesssim \e(t)t^{-|\alpha|},\quad |\alpha|\le 2\\
    |\partial_x^\alpha \partial_\xi^\beta a_0(t,x,\xi)|&\lesssim \e(t)t^{-1-\frac{|\alpha|}{2}},\quad
|\alpha|\ge 2.
  \end{split}
\end{equation}
This additional decay shall be used on time scales which are too small to allow $s^{\frac{1}{2}}\times s^{-\frac{1}{2}}$ packets at time $s$ to separate in time $t$.

Unfortunately, simply defining the parametrix $S_j^+$ by the evolution 
associated to the operator $|D|+a_{(0)}$ does not seem to work. Precisely,
the bounds \eqref{outl2j}, \eqref{finitel2j} and \eqref{freqloc}
appear to fail. This is because at each time $t$, there is leakage 
caused by the uncertainty principle to the regions appearing 
in \eqref{outl2j}, \eqref{finitel2j} and \eqref{freqloc}, which are 
outside the propagation region indicated by the Hamilton flow.
While this leakage does have rapid spatial decay, its time evolution
yields output which does not have the rapid decay in time as needed
in \eqref{outl2j}, \eqref{finitel2j} and \eqref{freqloc}.

Thus, in order to be able to prove the rapid $t$-decay in, e.g., \eqref{outl2j} \eqref{finitel2j} and \eqref{freqloc}, we shall introduce an artificial damping term $b_0\in l^1 S^{(1)}$, $b_0(t,x,\xi)\ge 0$. The role of $b_0$ is precisely 
to put a damping on the time evolution of the above mentioned leakeage.  
At the same time, $b_0$ is taken to be $0$ in the main propagation region.
We would like to define $S^+_j(t,s)$ to be the forward evolution operator
of the equation
\[ 
(D_t + |D| + a^w_0(t,x,D)) u = ib_0^w(t,x,D)u.
\]
However, in order to insure the frequency localization of our 
parametrix we replace $S^+_j(t,s)$ by the truncated operator
\[ 
S_{[-7,7]}(D_x)\cdot S_j^+.
\]
 We shall show that
\begin{align}
  \|x^\alpha S_{<-5}(D_x)S_j^+(t,s)P_j^+\|_{L^2\to L^2}&\lesssim (|t-s|+2^j)^{-N},\label{lowfreq}\\
  \|x^\alpha \partial^\beta S_{>5}(D_x)S_j^+(t,s)P_j^+\|_{L^2\to L^2}&\lesssim (|t-s|+2^j)^{-N},\label{highfreq}
\end{align}
and thus, the errors in \eqref{l2errorj} which result from this truncation are negligible.
We shall further prove the following bound on the
damping term
\begin{equation}
  \label{bbound}
  \begin{split}
    \|x^\alpha b_0^w(t,x,D)S_j^+(t,s)P_j^+\|_{L^2\to L^2}&\lesssim (|t-s|+2^j)^{-N},\\
    \|x^\alpha D_t b_0^w(t,x,D)S_j^+(t,s)P_j^+\|_{L^2\to L^2}&\lesssim (|t-s|+2^j)^{-N},
  \end{split}
\end{equation}
which shall yield \eqref{l2errorj}.

With $S_j^+(t,s)$ now fixed, property (iii) is trivial, and (i) follows from Proposition \ref{l2tol2}.
We proceed to the argument which yields our main 
pointwise bound \eqref{ptwisedecay}.  Here, we examine three cases separately.

{\bf Case 1:} $|t-s|\ge s$.  In this regime, we may neglect the damping.  For initial data
$u(s)=\delta_0$, we have
\[ T_{\frac{1}{s}}u(s,x_s,\xi)= s^{-\frac{n}{4}}e^{-\frac{x_s^2}{2s}}e^{ix_s\xi}. \]
Using Theorem \ref{tphase1}, we see that
\begin{align*}
|T_{\frac{1}{t}}u(t,x,\xi_t)|&\lesssim t^{-\frac{n}{4}}\int (1+t^{-1}|x-x_t(\xi_t,x_s)|^2)^{-N}
(1+s|\xi-\xi_s(x_s,\xi_t)|^2)^{-N} e^{-\frac{x_s^2}{2s}}\:dx_s\,d\xi \\
&\lesssim t^{-\frac{n}{4}}s^{-\frac{n}{2}} \int (1+t^{-1}|x-x_t(\xi_t,x_s)|^2)^{-N} e^{-\frac{x_s^2}{2s}}\:dx_s.
\end{align*}
For the remaining integral, we use that $x_s\to x_t(\xi_t,x_s)$ is Lipschitz.  See \eqref{xsxt1}.  Integrating
in $x_s$ then yields
\[
|T_{\frac{1}{t}}u(t,x,\xi_t)|\lesssim t^{-\frac{n}{4}}(1+t^{-1}|x-x_t(\xi_t,0)|^2)^{-N},\]
and by applying $T_{\frac{1}{t}}^*$, we have
\begin{align*}
|u(t,y)|&\lesssim t^{-\frac{n}{2}}\int (1+t^{-1}|x-x_t(\xi_t,0)|^2)^{-N} e^{-\frac{|y-x|^2}{2t}}\:dx\,d\xi_t\\
&\lesssim \int (1+t^{-1}|y-x_t(\xi_t,0)|^2)^{-N}\:d\xi_t.
\end{align*}
If $|t-s|\ge s$, then the map $\xi_t\to x_t(\xi_t,0)$ is zero homogeneous, by \eqref{xsxt1} 
has Lipschitz constant which is bounded by $t$, and has maximal rank $n-1$.  Hence, integration with respect
to $\xi_t$ yields
\[ |u(t,y)|\lesssim t^{-\frac{n-1}{2}}. \]

{\bf Case 2:} $1\le |t-s|\le s.$  Here, we reinitialize the time scale to prevent difficulties
which result from $s^{\frac{1}{2}}\times s^{-\frac{1}{2}}$ packets at time $s$ not separating
before time $t$.  In addition to \eqref{extraa}, we similarly require
\begin{equation}\label{extrab}
  \begin{split}
    |\partial_x^\alpha \partial_\xi^\beta b_0(t,x,\xi)|&\lesssim t^{-\frac{1}{2}-|\alpha|},
\quad |\alpha|\le 1,\\
|\partial_x^\alpha \partial_\xi^\beta b_0(t,x,\xi)|&\lesssim t^{-1-\frac{|\alpha|}{2}},\quad |\alpha|\ge 1
  \end{split}
\end{equation}
for $|t-s|<2^j$.  The additional regularity \eqref{extraa} and \eqref{extrab} is sufficient 
to show that $a_0$, $b_0$ remain in the appropriate symbol classes after the time translation which
sets the initial time to $t-s$.  Theorem \ref{tphase1} thus remains valid, and the bound follows
from the computation above in the translated coordinates.

{\bf Case 3:} $0\le |t-s|\le 1$.  Here, since our initial data is localized at frequency $1$, we may
simply use Sobolev embeddings combined with the $L^2$ bounds from Proposition \ref{l2tol2}:
\[
\|S^+_j(t,s)P_j^+ u_0\|_{L^\infty}\lesssim \| S^+_j(t,s)
P_j^+ u_0\|_{L^2}\lesssim \|P_j^+ u_0\|_{L^2} \lesssim \|u_0\|_{L^1}.\]

The rest of the proof is based on properties of $b_0$.  In particular, we use a construction which 
is quite similar to that of \cite{gS} to build a $b_0$ which allows us to prove the remaining required
estimates:
\eqref{outl2j}, \eqref{finitel2j}, \eqref{lowfreq}, \eqref{highfreq}, and \eqref{bbound}.  In particular, we have
\begin{lemma}\label{blemma}
There exists a symbol $b\in l^1 S^{(1)}$ which satisfies, in addition to \eqref{extrab},
\begin{enumerate}
  \item[(b1)] $t^{\frac{3}{4}}b$ is nonincreasing along the Hamiltonian flow for $D_t+|D_x|+a^w_0$, and
\[ 0<t^{\frac{3}{4}}b(t,x_t,\xi_t)<1 \implies b(2t,x_{2t},\xi_{2t})=0. \]
\item[(b2)] At the initial time, we have
\[ b(2^j, x,\xi)=0,\quad\text{ in }\quad \{2^{-3}<|\xi|<2^3,\quad 2^{j-2}<|x|<2^{j+2},\quad
x\xi > -2^{-4}|x|\}. \]
\item[(b3)] At any time $t\ge 2^j$, we have
\[ b(t,x,\xi)=t^{-\frac{3}{4}},\quad \text{outside }\quad \{2^{-4}<|\xi|<2^4\}\cap \{2^{-6}t<|x|<2^6t\}.\]
\end{enumerate}
\end{lemma}

Before proving this lemma, let us explain how such a damping term can
be used to complete the proof of Proposition \ref{K0}.  Indeed, we
have the following lemma which is essentially from \cite{gS}:
\begin{lemma}
  Assume that the symbol $b_0\in l^1S^{(1)}$ satisfies the properties
  (b1), (b2), and (b3) above.  Then, the bounds \eqref{outl2j},
  \eqref{finitel2j}, \eqref{freqloc}, \eqref{lowfreq}, \eqref{highfreq}, and
  \eqref{bbound} hold.
\end{lemma}

Indeed, once Theorem \ref{tphase1} has been established, the necessary
modifications to the arguments from \cite{gS} are quite simple.  It
only remains to establish the second estimate from \eqref{bbound}
which did not appear in \cite{gS}.  Here, however, we note
that, modulo negligible errors due to the frequency truncation of $S^+$,
\begin{multline*}
D_t b_0^w(t,x,D)S^+(t,s) = -i(\partial_t b_0)^w(t,x,D)S^+(t,s)
- b_0^w(t,x,D)|D_x|S^+(t,s)\\ - b_0^w(t,x,D)a^w_0(t,x,D)S^+(t,s)
+i b_0^w(t,x,D)b^w_0(t,x,D)S^+(t,s).\end{multline*}
Since the symbols of $(\partial_t b_0)^w$, $b_0^w |D_x|$, $b_0^w a_0^w$,
and $b_0^w b_0^w$ are all in $S^{(1)}_t$ and have supports which are contained in the support of $b_0$,
we may similarly apply Proposition 17 of \cite{gS} to obtain the estimate.

It now only remains to complete the construction of said damping terms $b$.
\begin{proof}[Proof of Lemma \ref{blemma}]
We define the increasing bounded function $e(s)$ by
\[
e(s) = \e^{-1} \int_{0}^s \frac{\e(\sigma)}{\sigma} d \sigma.
\]
%
Letting $\phi$ be a smooth, nondecreasing 
cutoff function which equals $0$ in $(-\infty,0)$ and $1$ in $(1,\infty)$, we
set
\[ b(t,x,\xi)=t^{-\frac{3}{4}}(1-\phi(b_1)\phi(b_2)\phi(b_3)\phi(b_4)\phi(b_5)) \]
with
\begin{itemize}
  \item {\em Cutoff frequencies which are too large}
\[ b_1(t,\xi)=\frac{2^{7/2}+e(t)-|\xi|}{\e(t)}, \]
  \item {\em Cutoff frequencies which are too small}
\[ b_2(t,\xi) = \frac{|\xi|-2^{-7/2}+c e(t)}{\e(t)} \]
where $c$ is a fixed small constant,
\item {\em Select outgoing waves}
\[ b_3(t,x,\xi) = \frac{2^{-\frac{1}2}|x||\xi|+ x\xi}{2^{-12}|x|},\]
\item {\em Cutoff values of $|x|$ which are too large}
\[ b_4(t,x)=\frac{2^6 t-|x|}{t}, \]
\item {\em Cutoff values of $|x|$ which are too small}
\[ b_5(t,x,\xi)=\frac{|x||\xi|-2^{-5}t|\xi|+x\xi}{2^{-10}t}.\]
\end{itemize}

We note that
\[ \{2^{-3}<|\xi|<2^3\}\cap \{2^{-2}t<|x|<2^2 t\}\cap \{x\xi>-2^{-4}|x|\}=D_t\subset \{b=0\} \]
if $\e$ is sufficiently small, while
\[ \{t^{\frac{3}{4}}b<1\}\subset E_t=\{2^{-4}<|\xi|<2^4\}\cap \{2^{-6}t<|x|<2^6t\}\cap \{x\xi >-2^{-1/2}|x||\xi|\}.\]
So, the conditions (b2) and (b3) are easily satisfied.

To prove (b1), it suffices to study the behavior of $b$ along the Hamilton flow within $E_t$ and
show that for each $b_j$, we have
\begin{equation}\label{bj}
\frac{d}{dt}b_j(t,x_t,\xi_t)\ge \frac{2}{t},\quad\text{in } E_t\cap \{0\le b_j\le 1\}.
\end{equation}
Here $t\to (x_t,\xi_t)$ now denotes a trajectory of the flow for $D_t+|D_x|+a_0^w$.
For $(x_t,\xi_t)\in E_t$, we have
\[\frac{d}{dt}\xi_t = O\Bigl(\frac{\e(t)}{t}\Bigr),\quad \frac{d}{dt}x_t = \frac{\xi_t}{|\xi_t|}+O(\e(t)).\]

We simply calculate
\[
\frac{d}{dt}  b_1(t,\xi_t)
\ge \frac{1}{\e t} - \frac{1}{t} - \frac{\e'(t)}{\e^2(t)}(2^{7/2}+e(t)-|\xi_t|) \ge \frac{2}{t},
\quad \text{in } \{0\le b_1\le 1\}
 \]
for $\e$ sufficiently small.  The computation for $b_2$ is identical.  For $b_3$, we have
\[ \frac{d}{dt}b_3(t,x_t,\xi_t)=\frac{|\xi_t|^2|x_t|^2 - (x_t\xi_t)^2}{2^{-12}|x_t|^3 |\xi_t|} + 
\frac{O(\e(t))}{t}\ge \frac{2}{t},\quad \text{in } E_t\cap \{0\le b_3\le 1\}.\]
For $b_4$, we compute
\[ \frac{d}{dt}b_4(t,x_t)=\frac{|x_t|^2|\xi_t| - tx_t\xi_t}{t^2|x_t||\xi_t|} + \frac{O(\e(t))}{t}\ge
\frac{2^5}{t},\quad\text{in } 0\le b_4\le 1. \]
Finally, for $b_5$ we also compute
\begin{align*}
  \frac{d}{dt}b_5(t,x_t,\xi_t)&=\frac{|x_t|^{-1}x_t\xi_t + |\xi_t|}{2^{-10}t}-\frac{|x_t||\xi_t|+x_t\xi_t}{2^{-10}
t^2} + \frac{O(\e(t))}{t}\\
&\ge \frac{2^5|\xi_t|}{t}\ge \frac{2}{t},\quad \text{in } E_t\cap \{0\le b_5\le 1\}.
\end{align*}

It remains to verify \eqref{extrab}, and hence that $b\in l^1S^{(1)}$, but this is straightforward.
\end{proof}

\bibliography{nls}

\begin{thebibliography}{10}

\bibitem{Alinhac}
Serge Alinhac.
\newblock On the {M}orawetz--{K}eel-{S}mith-{S}ogge inequality for the wave
  equation on a curved background.
\newblock {\em Publ. Res. Inst. Math. Sci.}, 42(3):705--720, 2006.

\bibitem{BT2}
Jean-Marc Bouclet and Nikolay Tzvetkov.
\newblock On global {S}trichartz estimates for non-trapping metrics.
\newblock {\em J. Funct. Anal.}, 254(6):1661--1682, 2008.

\bibitem{Brenner}
Philip Brenner.
\newblock On {$L\sb{p}-L\sb{p\sp{\prime} }$} estimates for the wave-equation.
\newblock {\em Math. Z.}, 145(3):251--254, 1975.

\bibitem{CoSaut}
Peter Constantin and Jean-Claude Saut.
\newblock Effets r\'egularisants locaux pour des \'equations dispersives
  g\'en\'erales.
\newblock {\em C. R. Acad. Sci. Paris S\'er. I Math.}, 304(14):407--410, 1987.

\bibitem{CKS}
Walter Craig, Thomas Kappeler, and Walter Strauss.
\newblock Microlocal dispersive smoothing for the {S}chr\"odinger equation.
\newblock {\em Comm. Pure Appl. Math.}, 48(8):769--860, 1995.

\bibitem{MR93i:35010}
Jean-Marc Delort.
\newblock {\em F.{B}.{I}. transformation}.
\newblock Springer-Verlag, Berlin, 1992.
\newblock Second microlocalization and semilinear caustics.

\bibitem{MR1795567}
Shin-ichi Doi.
\newblock Smoothing effects for {S}chr\"odinger evolution equation and global
  behavior of geodesic flow.
\newblock {\em Math. Ann.}, 318(2):355--389, 2000.

\bibitem{MR85f:35001}
Charles~L. Fefferman.
\newblock The uncertainty principle.
\newblock {\em Bull. Amer. Math. Soc. (N.S.)}, 9(2):129--206, 1983.

\bibitem{MR92k:22017}
Gerald~B. Folland.
\newblock {\em Harmonic analysis in phase space}.
\newblock Princeton University Press, Princeton, NJ, 1989.

\bibitem{MR1151250}
J.~Ginibre and G.~Velo.
\newblock Smoothing properties and retarded estimates for some dispersive
  evolution equations.
\newblock {\em Comm. Math. Phys.}, 144(1):163--188, 1992.

\bibitem{MR2131050}
Andrew Hassell, Terence Tao, and Jared Wunsch.
\newblock A {S}trichartz inequality for the {S}chr\"odinger equation on
  nontrapping asymptotically conic manifolds.
\newblock {\em Comm. Partial Differential Equations}, 30(1-3):157--205, 2005.

\bibitem{KSS}
Markus Keel, Hart~F. Smith, and Christopher~D. Sogge.
\newblock Almost global existence for some semilinear wave equations.
\newblock {\em J. Anal. Math.}, 87:265--279, 2002.
\newblock Dedicated to the memory of Thomas H.\ Wolff.

\bibitem{MR1646048}
Markus Keel and Terence Tao.
\newblock Endpoint {S}trichartz estimates.
\newblock {\em Amer. J. Math.}, 120(5):955--980, 1998.

\bibitem{KPV}
Carlos~E. Kenig, Gustavo Ponce, and Luis Vega.
\newblock On the {Z}akharov and {Z}akharov-{S}chulman systems.
\newblock {\em J. Funct. Anal.}, 127(1):204--234, 1995.

\bibitem{MR2094851}
Herbert Koch and Daniel Tataru.
\newblock Dispersive estimates for principally normal pseudodifferential
  operators.
\newblock {\em Comm. Pure Appl. Math.}, 58(2):217--284, 2005.

\bibitem{laxphil}
Peter~D. Lax and Ralph~S. Phillips.
\newblock {\em Scattering theory}, volume~26 of {\em Pure and Applied
  Mathematics}.
\newblock Academic Press Inc., Boston, MA, second edition, 1989.
\newblock With appendices by Cathleen S. Morawetz and Georg Schmidt.

\bibitem{MMT}
Jeremy Marzuola, Jason Metcalfe, and Daniel Tataru.
\newblock Strichartz estimates and local smoothing estimates for
  asympototically flat {S}chr\"odinger equations.
\newblock {\em J. Funct. Anal.}, 255(6):1497--1553, 2008.

\bibitem{MetSo}
Jason Metcalfe and Christopher~D. Sogge.
\newblock Long-time existence of quasilinear wave equations exterior to
  star-shaped obstacles via energy methods.
\newblock {\em SIAM J. Math. Anal.}, 38(1):188--209 (electronic), 2006.

\bibitem{MR1168960}
Gerd Mockenhaupt, Andreas Seeger, and Christopher~D. Sogge.
\newblock Local smoothing of {F}ourier integral operators and
  {C}arleson-{S}j\"olin estimates.
\newblock {\em J. Amer. Math. Soc.}, 6(1):65--130, 1993.

\bibitem{morawetz}
Cathleen~S. Morawetz.
\newblock Time decay for the nonlinear {K}lein-{G}ordon equations.
\newblock {\em Proc. Roy. Soc. Ser. A}, 306:291--296, 1968.

\bibitem{RZ}
Luc Robbiano and Claude Zuily.
\newblock Strichartz estimates for {S}chr\"odinger equations with variable
  coefficients.
\newblock {\em M\'em. Soc. Math. Fr. (N.S.)}, (101-102):vi+208, 2005.

\bibitem{RT}
I.~Rodnianski and T.~Tao.
\newblock Longtime decay estimates for the {S}chr\"odinger equation on
  manifolds.
\newblock In {\em Mathematical aspects of nonlinear dispersive equations},
  volume 163 of {\em Ann. of Math. Stud.}, pages 223--253. Princeton Univ.
  Press, Princeton, NJ, 2007.

\bibitem{Sjo}
Per Sj{\"o}lin.
\newblock Regularity of solutions to the {S}chr\"odinger equation.
\newblock {\em Duke Math. J.}, 55(3):699--715, 1987.

\bibitem{Sj1}
Johannes Sj{\"o}strand.
\newblock Singularit\'es analytiques microlocales.
\newblock In {\em Ast\'erisque, 95}, pages 1--166. Soc. Math. France, Paris,
  1982.

\bibitem{MR1644105}
Hart~F. Smith.
\newblock A parametrix construction for wave equations with {$C\sp {1,1}$}
  coefficients.
\newblock {\em Ann. Inst. Fourier (Grenoble)}, 48(3):797--835, 1998.

\bibitem{SS}
Hart~F. Smith and Christopher~D. Sogge.
\newblock On {S}trichartz and eigenfunction estimates for low regularity
  metrics.
\newblock {\em Math. Res. Lett.}, 1(6):729--737, 1994.

\bibitem{smithsogge}
Hart~F. Smith and Christopher~D. Sogge.
\newblock Global {S}trichartz estimates for nontrapping perturbations of the
  {L}aplacian.
\newblock {\em Comm. Partial Differential Equations}, 25(11-12):2171--2183,
  2000.

\bibitem{MR1909638}
Hart~F. Smith and Daniel Tataru.
\newblock Sharp counterexamples for {S}trichartz estimates for low regularity
  metrics.
\newblock {\em Math. Res. Lett.}, 9(2-3):199--204, 2002.

\bibitem{MR1924470}
Gigliola Staffilani and Daniel Tataru.
\newblock Strichartz estimates for a {S}chr\"odinger operator with nonsmooth
  coefficients.
\newblock {\em Comm. Partial Differential Equations}, 27(7-8):1337--1372, 2002.

\bibitem{Strauss}
Walter~A. Strauss.
\newblock Dispersal of waves vanishing on the boundary of an exterior domain.
\newblock {\em Comm. Pure Appl. Math.}, 28:265--278, 1975.

\bibitem{Str}
Robert~S. Strichartz.
\newblock Restrictions of {F}ourier transforms to quadratic surfaces and decay
  of solutions of wave equations.
\newblock {\em Duke Math. J.}, 44(3):705--714, 1977.

\bibitem{nlw}
Daniel Tataru.
\newblock Strichartz estimates for operators with nonsmooth coefficients and
  the nonlinear wave equation.
\newblock {\em Amer. J. Math.}, 122(2):349--376, 2000.

\bibitem{cs}
Daniel Tataru.
\newblock Strichartz estimates for second order hyperbolic operators with
  nonsmooth coefficients. {I}{I}.
\newblock {\em Amer. J. Math.}, 123(3):385--423, 2001.

\bibitem{MR1944027}
Daniel Tataru.
\newblock On the {F}efferman-{P}hong inequality and related problems.
\newblock {\em Comm. Partial Differential Equations}, 27(11-12):2101--2138,
  2002.

\bibitem{lp}
Daniel Tataru.
\newblock Strichartz estimates for second order hyperbolic operators with
  nonsmooth coefficients. {III}.
\newblock {\em J. Amer. Math. Soc.}, 15(2):419--442 (electronic), 2002.

\bibitem{phasespace}
Daniel Tataru.
\newblock Phase space transforms and microlocal analysis.
\newblock In {\em Phase space analysis of partial differential equations. Vol.
  II}, Pubbl. Cent. Ric. Mat. Ennio Giorgi, pages 505--524. Scuola Norm. Sup.,
  Pisa, 2004.

\bibitem{gS}
Daniel Tataru.
\newblock Parametrices and dispersive estimates for {S}chr\"odinger operators
  with variable coefficients.
\newblock {\em Amer. J. Math.}, 130(3):571--634, 2008.

\bibitem{MR1766415}
Michael~E. Taylor.
\newblock {\em Tools for {PDE}}, volume~81 of {\em Mathematical Surveys and
  Monographs}.
\newblock American Mathematical Society, Providence, RI, 2000.
\newblock Pseudodifferential operators, paradifferential operators, and layer
  potentials.

\bibitem{Veg}
Luis Vega.
\newblock Schr\"odinger equations: pointwise convergence to the initial data.
\newblock {\em Proc. Amer. Math. Soc.}, 102(4):874--878, 1988.

\end{thebibliography}

\end{document}